\documentclass{article}
\usepackage{makeidx}
\usepackage{amssymb}
\usepackage{graphicx}
\usepackage{amsmath}

\newtheorem{theorem}{Theorem}

\newtheorem{corollary}[theorem]{Corollary}

\newtheorem{definition}[theorem]{Definition}

\newtheorem{lemma}[theorem]{Lemma}

\newtheorem{proposition}[theorem]{Proposition}

\newtheorem{remark}[theorem]{Remark}
\newtheorem{rem}[theorem]{Remark}

\input{tcilatex}
\makeindex

\begin{document}

\title{Solitons in Schr\"{o}dinger-Maxwell equations}
\author{Vieri Benci\thanks{
Dipartimento di Matematica, Universit\`a degli Studi di Pisa, Via F.
Buonarroti 1/c, Pisa, ITALY and Department of Mathematics, College of
Science, King Saud University, Riyadh, 11451, SAUDI ARABIA. e-mail: \texttt{%
benci@dma.unipi.it}} \and Donato Fortunato\thanks{%
Dipartimento di Matematica, Universit\`{a} degli Studi di Bari "Aldo Moro"
and INFN sezione di Bari, Via Orabona 4, Bari, ITALY. e-mail: \texttt{\
fortunat@dm.uniba.it, }}}
\maketitle

\begin{abstract}
In this paper we study the Nonlinear Schr\"{o}dinger-Maxwell equations
(NSM). We are interested to analyse the existence of solitons, namely of
finite energy solutions which exhibit stability properties. This paper is
divided in two parts. In the first, we give an abstract definition of
soliton and we develope an abstract existence theory. In the second, we
apply this theory to NSM.

\ 

\textbf{Key words:} Maxwell equations, Nonlinear Schr\"{o}dinger equation,
solitary waves, Hylomorphic solitons, variational methods.

\ 

\textbf{AMS Subject Classification}: 35Q61, 35Q55, 37K40, 35A15.
\end{abstract}

\tableofcontents

\section{Introduction}

In this paper we study a system of equations obtained by coupling the Schr%
\"{o}dinger equation with the Maxwell equations (NSM) (see eq. (\ref{aaa}), (%
\ref{bbb})). This system, usually called the Schr\"{o}dinger-Poisson system
or Schr\"{o}dinger Maxwell system, describes many interesting physical
situations (see e.g. \cite{sanchez} and its references). We are interested
to analyse the existence of solitons, namely of finite energy solutions
which exhibit a strong form of stability. In particular we are interested in
a class of solitons which , following \cite{BBBM}, \cite{hylo}, \cite{milan}%
, \cite{bfniren}, are called hylomorphic. The existence of such solitons is
due to the interplay between two constants of the motion: the energy and the
charge.

This paper is divided in two parts.

In the first part, following \cite{befocao}, we give an abstract definition
of soliton and we develope an abstract existence theory for hylomorphic
solitons. This theory is based on concentration-compactness type arguments
(see \cite{Li84a}, \cite{Li84b}).

In the second part this theory has been used to prove the existence of
hylomorphic solitons for NSM (see Theorems \ref{exnsm} and \ref{nsmbis})
when the coupling constant $q$ is sufficiently small. If $q$ $=0$ the NSM
reduce to the Schr\"{o}dinger equation. So Theorems \ref{exnsm} and \ref%
{nsmbis} extend to the case of NSM some of the well known stability results
stated for the Schr\"{o}dinger equation (see e. g. \cite{Cazli}, \cite{We86}%
, \cite{GSS90}, \cite{gss87}, \cite{BBGM}, \cite{sulesulem} and its
references). NSM has been largely studied by many authors and under various
assumptions on the nonlinear term. There is a huge bibliography on this
subject and the list of our references is far to be complete. For the
existence of solutions we refer to \cite{ambr}, \cite{ambruiz}, \cite{belsic}%
, \cite{befosch}, \cite{candela}, \cite{coclit}, \cite{coc}, \cite%
{ceramvaira}, \cite{tea}, \cite{davpompva}, \cite{dav2}, \cite{kikuki}, \cite%
{pisanisic}, \cite{salvatore}, \cite{ruiz1}, \cite{sanchez}. However we know
only few results (\cite{belsic}, \cite{kikukibis}) proving the existence of
stable solitary waves (namely solitons) for such equations. For the study of
some qualitative properties of the solutions, like the presence of
concentration phenomena and the study of semiclassical limits, we refer to 
\cite{teawei}, \cite{teawei2}, \cite{tea1s}, \cite{iannivair}, \cite{Ruiz}, 
\cite{ruizvaira}.

\bigskip

Our approach to NSM presents the following novelties:

\begin{itemize}
\item The proof of the existence result is based on a new abstract framework.

\item The nonlinear term is not assumed to be homogeneous.

\item The stability of the solutions is proved.

\item The presence of a "lattice type" potential $V(x)$ is allowed.
\end{itemize}

\section{Solitary waves and solitons: abstract theory\label{ab}}

In this section, following \cite{befocao} and \cite{befolib}, we introduce a
functional abstract framework which allows to define solitary waves,
solitons and hylomorphic solitons. Then, we will state some abstract
existence theorems. These theorems are based on a general minimization
principle related to the concentration compactness techniques.

\subsection{Basic definitions}

\textit{Solitary waves} and solitons are particular \textit{states} of a
dynamical system described by one or more partial differential equations.
Thus, we assume that the states of this system are described by one or more 
\textit{fields} which mathematically are represented by functions 
\begin{equation*}
\mathbf{u}:\mathbb{R}^{N}\rightarrow V
\end{equation*}%
where $V$ is a vector space with norm $\left\vert \ \cdot \ \right\vert _{V}$
which is called the internal parameters space. We assume the system to be
deterministic; this means that it can be described as a dynamical system $%
\left( X,\gamma \right) $ where $X$ is the set of the states and $\gamma :%
\mathbb{R}\times X\rightarrow X$ is the time evolution map. If $\mathbf{u}%
_{0}(x)\in X,$ the evolution of the system will be described by the function 
\begin{equation}
\mathbf{u}\left( t,x\right) :=\gamma _{t}\mathbf{u}_{0}(x).  \label{flusso}
\end{equation}%
We assume that the states of $X$ have "finite energy" so that they decay at $%
\infty $ sufficiently fast and that%
\begin{equation}
X\subset L_{loc}^{1}\left( \mathbb{R}^{N},V\right) .  \label{lilla}
\end{equation}

Using this framework, we give the following definitions:\label{pag}

\begin{definition}
A dynamical system $\left( X,\gamma \right) $ is called of FT type
(field-theory-type) if $X$ is a Hilbert space of functions satisfying (\ref%
{lilla}).
\end{definition}

For every $\tau \in \mathbb{R}^{N},$ and $\mathbf{u}\in X$, we set 
\begin{equation}
\left( T_{\tau }\mathbf{u}\right) \left( x\right) =\mathbf{u}\left( x-\tau
\right) .  \label{ggg}
\end{equation}%
Clearly, the group 
\begin{equation}
\mathcal{T}=\left\{ T_{\tau }|\ \tau \in \mathbb{R}^{N}\right\} ;  \label{gg}
\end{equation}%
is a representation of the group of translations%
\index{translations}.

\begin{definition}
A set $\Gamma \subset X$ is called compact up to space tanslations or $%
\mathcal{T}$-compact%
\index{T-compact} if for any sequence $\mathbf{u}_{n}(x)\in \Gamma \ $there
is a subsequence $\mathbf{u}_{n_{k}}$ and a sequence $\tau _{k}\in \mathbb{R}%
^{N}$ such that $\mathbf{u}_{n_{k}}(x-\tau _{k})$ is convergent.
\end{definition}

Now, we want to give an abstract definition of solitary wave. Roughly
speaking a solitary wave is a field whose energy travels as a localized
packet and which preserves this localization in time. For example, consider
a solution of a field equation having the following form: 
\begin{equation}
\mathbf{u}\left( t,x\right) =u(x-vt-x_{0})e^{i(v\mathbf{\cdot }x\mathbf{-}%
\omega t)};\ u\in L^{2}(\mathbb{R}^{N}).  \label{solwav}
\end{equation}%
The field (\ref{solwav}) is a solitary wave depending on the constants $%
x_{0},v$ and $\omega $. The evolution of a solitary wave is a translation
plus a mild change of the internal parameters (in this case the phase).

This situation can be formalized by the following definition:

\begin{definition}
\label{solw} If $\mathbf{u}\in X,$ we denote the closure of the orbit of $%
\mathbf{u}$ by%
\begin{equation*}
\mathcal{O}\left( \mathbf{u}\right) :=%
\overline{\left\{ \gamma _{t}\mathbf{u}(x)\ |\ t\in \mathbb{R}\right\} }.
\end{equation*}%
A state $\mathbf{u}\in X$ is called solitary wave if

\begin{itemize}
\item (i) $0\notin \mathcal{O}\left( \mathbf{u}\right) ;$

\item (ii) $\mathcal{O}\left( \mathbf{u}\right) $ is $\mathcal{T}$-compact%
\index{T-compact}.
\end{itemize}
\end{definition}

Clearly, (\ref{solwav}) describes a solitary wave according to the
definition above. The standing waves, namely objects of the form 
\begin{equation}
\gamma _{t}\mathbf{u=u}(t,x)=u(x)e^{-i\omega t},\ \ u\in L^{2}(\mathbb{R}%
^{N}),\ u\neq 0  \label{stand}
\end{equation}%
probably are the "simplest" solitary waves. In this case the orbit $\mathcal{%
O}\left( \mathbf{u}\right) $ is compact.

Take $X=L^{1}(\mathbb{R}^{N})$\ and\ $\mathbf{u}\in X;$ if $\gamma _{t}%
\mathbf{u}=\mathbf{u}\left( e^{t}x\right) $, $\mathbf{u}$ is not a solitary
wave since $\left\Vert \gamma _{t}\mathbf{u}\right\Vert _{X}\mathbf{%
\rightarrow }0$ as $t\rightarrow +\infty $ and (i) is clearly violated. If $%
\gamma _{t}\mathbf{u}=e^{t}\mathbf{u}\left( e^{t}x\right) ,$ $\mathbf{u}$ is
not a solitary wave since (ii) of Def. \ref{solw} does not hold. Also,
according to our definition, a "couple" of solitary waves is not a solitary
wave: in fact 
\begin{equation*}
\gamma _{t}\mathbf{u}=\left[ u(x-vt)+u(x+vt)\right] e^{i(v\mathbf{\cdot }x%
\mathbf{-}\omega t)},\ u\in L^{2}(\mathbb{R}^{N})
\end{equation*}%
is not a solitary wave since (ii) is violated.

The \textit{solitons%
\index{soliton}} are solitary waves characterized by some form of stability.
To define them at this level of abstractness, we need to recall some well
known notions in the theory of dynamical systems.

\begin{definition}
A set $\Gamma \subset X$ is called \textit{invariant%
\index{invariant}} if $\forall \mathbf{u}\in \Gamma ,\forall t\in \mathbb{R}%
,\ \gamma _{t}\mathbf{u}\in \Gamma .$
\end{definition}

\begin{definition}
Let $\left( X,d\right) $ be a metric space and let $\left( X,\gamma \right) $
be a dynamical system. An invariant set $\Gamma \subset X$ is called stable%
\index{stable invariant set}, if $\forall \varepsilon >0,$ $\exists \delta
>0,\;\forall \mathbf{u}\in X$, 
\begin{equation*}
d(\mathbf{u},\Gamma )\leq \delta ,
\end{equation*}%
implies that 
\begin{equation*}
\forall t\geq 0,%
\text{ }d(\gamma _{t}\mathbf{u,}\Gamma )\leq \varepsilon .
\end{equation*}
\end{definition}

Now we are ready to give the definition of soliton:

\begin{definition}
\label{ds} A state $\mathbf{u}\in X$ is called soliton if $\mathbf{u}\in
\Gamma \subset X$ where

\begin{itemize}
\item (i) $\Gamma $ is an invariant, stable set

\item (ii) $\Gamma $ is $\mathcal{T}$-compact

\item (iii) $0\notin \Gamma $.
\end{itemize}

The set $\Gamma $ is called soliton manifold.
\end{definition}

The above definition needs some explanation. First of all notice that every $%
\mathbf{u}\in \Gamma $ is a soliton and that every soliton is a solitary
wave. Now for simplicity, we assume that $\Gamma $ is a manifold\footnote{%
actually, in many concrete models, this is the generic case; this is the
reason why $\Gamma $ is called \textit{soliton manifold} even if it might
happen that it is not a manifold.}. Then (ii) implies that $\Gamma $ is
finite dimensional. Since $\Gamma $ is invariant, $\mathbf{u}\in \Gamma
\Rightarrow \gamma _{t}\mathbf{u}\in \Gamma $ for every time. Thus, since $%
\Gamma $ is finite dimensional, the evolution of $\mathbf{u}$ is described
by a finite number of parameters$.$ The dynamical system $\left( \Gamma
,\gamma \right) $\ behaves as a point in a finite dimensional phase space.
By the stability of $\Gamma $, a small perturbation of $\mathbf{u}$ remains
close to $\Gamma .$ However, in this case, its evolution depends on an
infinite number of parameters. Thus, this system appears as a finite
dimensional system with a small perturbation.

We now assume that the dynamical system $\left( X,\gamma \right) $ has two
constants of motion: the energy $E$ and the hylenic charge $C.$ At this
level of abstraction, the name energy and hylenic charge are conventional,
but in the applications,\ $E$ and $C$ will be the energy and the hylenic
charge as defined in section \ref{ec}.

\begin{definition}
\label{tdc}A state $\mathbf{u}_{0}\in X$ is called \textbf{hylomorphic
soliton%
\index{soliton!hylomorphic}} if it is a soliton according to Def. \ref{ds}
and if the soliton manifold $\Gamma $ has the following structure%
\begin{equation}
\Gamma =\Gamma \left( e_{0},c_{0}\right) =\left\{ \mathbf{u}\in X\ |\ E(%
\mathbf{u})=e_{0},\ \left\vert C(\mathbf{u})\right\vert =c_{0}\right\}
\label{plis}
\end{equation}%
where%
\begin{equation}
e_{0}=\min \left\{ E(\mathbf{u})\ |\ \left\vert C(\mathbf{u})\right\vert
=c_{0}\right\} .  \label{minbis}
\end{equation}
\end{definition}

Notice that, by (\ref{minbis}), we have that a hylomorphic soliton $\mathbf{u%
}_{0}$ satisfies the following nonlinear eigenvalue problem:%
\begin{equation*}
E^{\prime }(\mathbf{u}_{0})=\lambda C^{\prime }(\mathbf{u}_{0}).
\end{equation*}%
In general, a minimizer $\mathbf{u}_{0}$ of $E$ on%
\begin{equation*}
\mathfrak{M}_{c_{_{0}}}:=\left\{ \mathbf{u}\in X\ |\ \left\vert C(\mathbf{u}%
)\right\vert =c_{0}\right\} ,
\end{equation*}%
is not a soliton; in fact, according to Def. \ref{ds}, it is necessary to
check the following facts:

\begin{itemize}
\item (i) the set $\Gamma \left( e_{0},c_{0}\right) $ is stable.

\item (ii) the set $\Gamma \left( e_{0},c_{0}\right) $ is $\mathcal{T}$%
-compact%
\index{T-compact} (i.e. compact up to translations).

\item (iii) $0\notin \Gamma \left( e_{0},c_{0}\right) $ since otherwise,
some $\mathbf{u}\in \Gamma \left( e_{0},c_{0}\right) $ is not even a
solitary wave (see Def. \ref{solw} (i)).
\end{itemize}

In concrete cases, the point (i) is the most delicate point to prove. If (i)
does not hold, according to our definitions, $\mathbf{u}_{0}$ is a solitary
wave but not a soliton.

\subsection{ An abstract minimization theorem}

In the previous section, we have seen that the existence of hylomorphic
soliton is related to the existence of minimizers of the energy. So in this
section we assume that $X$ is a Hilbert space and that $E$ and $C$ are two
differentiable functionals defined on it and we will investigate the
following minimization problem%
\begin{equation}
\min_{\mathbf{u}\in \mathfrak{M}_{c}}E(\mathbf{u})\ \ 
\text{where\ \ }\mathfrak{M}_{c}:=\left\{ \mathbf{u}\in X\ |\ \left\vert C(%
\mathbf{u})\right\vert =c\right\} .  \label{MinPr}
\end{equation}

\subsubsection{Preliminary notions}

We need a few abstract definitions some of which have been introduced in 
\cite{befocao}. In the following $G$ will denote a group with a unitary
action on $X.$

\begin{definition}
A subset $\Gamma \subset X$ is called $G$-invariant%
\index{G-invariant set} if 
\begin{equation*}
\forall \mathbf{u}\in \Gamma ,\ \forall g\in G,\ g\mathbf{u}\in \Gamma .
\end{equation*}
\end{definition}

In many concrete situations, $G$ will be a subgroup of the translations
group $\mathcal{T}$.

\begin{definition}
A sequence $\mathbf{u}_{n}$ in $X$ is called $G$\emph{-compact }%
\index{G-compact sequence}if there is a subsequence $\mathbf{u}_{n_{k}}$ and
a sequence $g_{k}\in G$ such that $g_{k}\mathbf{u}_{n_{k}}$ is convergent. A
subset $\Gamma \subset X$ is called $G$-compact if every sequence in $\Gamma 
$ is $G$\emph{-compact.}
\end{definition}

If $G=\left\{ Id\right\} $ or more in general it is a compact group, $G$%
\emph{-}compactness implies compactness. If $G$ is not compact such as the
translation group $\mathcal{T}$, $G$\emph{-}compactness is a weaker notion
than compactness.

\begin{definition}
\label{gcompattoa}A $G$-invariant functional%
\index{G-compact functional} $J$ on $X$ is called $G$-compact if any
minimizing sequence $\mathbf{u}_{n}$ is $G$-compact.
\end{definition}

Clearly a $G$-compact functional has a $G$-compact set of minimizers.

\begin{definition}
\label{sp}We say that a functional $F$ on $X$ has the splitting property%
\index{splitting property} if given a sequence $\mathbf{u}_{n}=\mathbf{u}+%
\mathbf{w}_{n}\in X$ such that $\mathbf{w}_{n}$ converges weakly to $0$, we
have that 
\begin{equation*}
F(\mathbf{u}_{n})=F(\mathbf{u})+F(\mathbf{w}_{n})+o(1)
\end{equation*}
\end{definition}

\begin{remark}
\label{quadratic}Every continuous quadratic form satisfies the splitting
property; in fact, in this case, we have that $F(\mathbf{u}):=\left\langle L%
\mathbf{u},\mathbf{u}\right\rangle $ for some continuous selfajoint operator 
$L;$ then, given a sequence $\mathbf{u}_{n}=\mathbf{u}+\mathbf{w}_{n}$ with $%
\mathbf{w}_{n}\rightharpoonup 0$ weakly, we have that%
\begin{align*}
F(\mathbf{u}_{n})& =\left\langle L\mathbf{u},\mathbf{u}\right\rangle
+\left\langle L\mathbf{w}_{n},\mathbf{w}_{n}\right\rangle +2\left\langle L%
\mathbf{u},\mathbf{w}_{n}\right\rangle \\
& =F(\mathbf{u})+F(\mathbf{w}_{n})+o(1)
\end{align*}
\end{remark}

\begin{definition}
\label{na}\textbf{\ }A sequence $\mathbf{u}_{n}\in X$ is called \textbf{%
vanishing sequence}\textit{\ if it is bounded and if for any sequence }$%
g_{n}\in G$ the sequence $g_{n}\mathbf{u}_{n}$ converges weakly to $0.$
\end{definition}

So, if $\mathbf{u}_{n}\rightarrow 0$ strongly, $\mathbf{u}_{n}$ is a
vanishing sequence. However, if $\mathbf{u}_{n}\rightharpoonup 0$ weakly, it
might happen that it is not a vanishing sequence; namely it might exist a
subsequence $\mathbf{u}_{n_{k}}$ and a sequence $g_{k}\in G$ such that $g_{k}%
\mathbf{u}_{n_{k}}$ is weakly convergent to some $\mathbf{%
\bar{u}}\neq 0$. Let see an example; if $u_{0}\in X\subset L^{1}\left( 
\mathbb{R}^{N}\right) $ and $x_{n}\rightarrow +\infty ,$ then the sequence $%
T_{x_{n}}u_{0}=u_{0}(x-x_{n})$ is not vanishing. Clearly, in this example $G$
contains the group of translations (\ref{gg}).

\bigskip

Now, we set 
\begin{equation}
\Lambda \left( \mathbf{u}\right) :=\frac{E\left( \mathbf{u}\right) }{%
\left\vert C\left( \mathbf{u}\right) \right\vert },  \label{lambda}
\end{equation}%
$\Lambda $ will be called \textbf{hylenic ratio}.

The notions of vanishing sequence and of hylenic ratio allow to introduce
the following (important) definition:

\begin{definition}
\label{dhc}We say that the \textit{hylomorphy condition holds%
\index{hylomorphy condition} }if%
\begin{equation}
\underset{\mathbf{u}\in X}{\inf }%
\frac{E\left( \mathbf{u}\right) }{\left\vert C\left( \mathbf{u}\right)
\right\vert }<\Lambda _{0}.  \label{hh}
\end{equation}%
where%
\begin{equation}
\Lambda _{0}:=\ \inf \left\{ \lim \inf \ \Lambda (\mathbf{u}_{n})\ |\ 
\mathbf{u}_{n}\ \text{is a vanishing sequence}\right\}  \label{hylo}
\end{equation}%
Moreover, we say that $\mathbf{u}_{0}\in X$ satisfies the \textit{hylomorphy
condition }if, 
\begin{equation}
\frac{E\left( \mathbf{u}_{0}\right) }{\left\vert C\left( \mathbf{u}%
_{0}\right) \right\vert }<\Lambda _{0}.  \label{yc}
\end{equation}
\end{definition}

So, if $\mathbf{u}_{n}$ is a bounded sequence, we have the following:%
\begin{equation*}
\lim \inf \Lambda \left( \mathbf{u}_{n}\right) <\Lambda _{0}\Rightarrow
\exists \mathbf{u}_{n_{k}},g_{k}\in G:g_{k}\mathbf{u}_{n_{k}}\rightharpoonup 
\mathbf{\bar{u}}\neq 0.
\end{equation*}

In order to apply the existence theorems of the next subsection, it is
necessary to estimate $\Lambda _{0};$ the following propositons may help to
do this.

\begin{proposition}
\label{diesis}Assume that there exists a seminorm $\left\Vert \cdot
\right\Vert _{\sharp }$ on $X$ such that\label{kk} 
\begin{equation}
\left\{ \mathbf{u}_{n}\ \text{is a vanishing sequence}\right\} \Rightarrow
\left( \left\Vert \mathbf{u}_{n}\right\Vert _{\sharp }\rightarrow 0\right)
\label{seminorm}
\end{equation}%
Then%
\begin{equation}
\underset{\left\Vert \mathbf{u}\right\Vert _{\sharp }\rightarrow 0}{\lim
\inf }\ \Lambda (\mathbf{u})\leq \Lambda _{0}\leq \ \underset{\left\Vert 
\mathbf{u}\right\Vert \rightarrow 0}{\lim \inf }\ \Lambda (\mathbf{u}).
\label{aaaa}
\end{equation}
\end{proposition}

\textbf{Proof. }By definition \ref{na} and by (\ref{seminorm}) we have 
\begin{equation}
\left( \left\Vert \mathbf{u}_{n}\right\Vert \rightarrow 0\right) \Rightarrow
\left( \mathbf{u}_{n}\ \text{ vanishing sequence}\right) \Rightarrow \left(
\left\Vert \mathbf{u}_{n}\right\Vert _{\sharp }\rightarrow 0\right)
\label{hylo2}
\end{equation}

Then, by (\ref{hylo}) and (\ref{hylo2}), we get (\ref{aaaa}).

\bigskip $\square $

\subsubsection{ The minimization result \label{33}}

We shall make the following assumptions on the functionals $E$ and $C$:

\begin{itemize}
\item \textit{(EC-0) \textbf{(Values at 0)}}%
\begin{equation*}
E(0)=C(0)=0;\ E^{\prime }(0)=C^{\prime }(0)=0.\ 
\end{equation*}

\item \textit{(EC-1)\textbf{(Invariance)} }$E(\mathbf{u})$\textit{\ and }$C(%
\mathbf{u})$\textit{\ are }$G$\textit{-invariant.}

\item \textit{(EC-2)\textbf{(Splitting property)} }$E$\textit{\ and }$C$%
\textit{\ satisfy the splitting property (see Definition \ref{sp}).}

\item \textit{(EC-3)\textbf{(Coercivity)\label{coercivity} }We assume that} 
\textit{there exists }$a>0$ \textit{and }$s>1$\textit{\ such that}

\begin{itemize}
\item (i) $\forall \mathbf{u}\neq 0,\ C(\mathbf{u})>0\ $\textit{and} $E(%
\mathbf{u})+aC(\mathbf{u})^{s}>0;$

\item (ii) \textit{if }$\left\Vert \mathbf{u}\right\Vert \rightarrow \infty
,\ $\textit{then} $E(\mathbf{u})+aC(\mathbf{u})^{s}\rightarrow \infty ;$

\item (iii) \textit{for any} \textit{bounded sequence }$\mathbf{u}_{n}$ 
\textit{in }$X$ \textit{such that} $E(\mathbf{u}_{n})+aC(\mathbf{u}%
_{n})^{s}\rightarrow 0,\ $\textit{we have that }$\mathbf{u}_{n}\rightarrow
0. $
\end{itemize}
\end{itemize}

Now we can state the main results. We start with a technical lemma.

\begin{lemma}
\label{leman} Assume that (EC-2) and (EC-3) (i) are satisfied. Let $\mathbf{u%
}_{n}=\mathbf{u}+\mathbf{w}_{n}\in X$ be a sequence such that $\mathbf{u}%
\neq 0,\mathbf{w}_{n}\neq 0\mathbf{\ }$and $\mathbf{w}_{n}$ converges weakly
to $0$. Then, up to a subsequence, we have 
\begin{equation}
\lim \Lambda \left( \mathbf{u+\mathbf{w}_{n}}\right) \geq \min \left(
\Lambda \left( \mathbf{u}\right) ,\lim \Lambda \left( \mathbf{\mathbf{w}_{n}}%
\right) \right)  \label{tiga}
\end{equation}
\end{lemma}

\textbf{Proof.} The proof is contained in \cite{befocao}. We shall repeat it
for completeness. Given four real numbers $A,B,a,b,$ (with $B,b>0$), we have
that%
\begin{equation}
\frac{A+a}{B+b}\geq \min \left( \frac{A}{B},\frac{a}{b}\right)  \label{figa}
\end{equation}%
In fact, suppose that $\frac{A}{B}\geq \frac{a}{b};$ then%
\begin{equation*}
\frac{A+a}{B+b}=\frac{\frac{A}{B}B+\frac{a}{b}b}{B+b}\geq \frac{\frac{a}{b}B+%
\frac{a}{b}b}{B+b}=\frac{a}{b}\geq \min \left( \frac{A}{B},\frac{a}{b}\right)
\end{equation*}%
Notice that the equality holds if and only if 
\begin{equation}
\frac{A}{B}=\frac{a}{b}.  \label{abAB}
\end{equation}%
Since $\mathbf{u}\neq 0\ $and $\mathbf{w}_{n}\neq 0\mathbf{,\ }$by (EC-3)
(i), we have $C\left( \mathbf{u)+}C(\mathbf{\mathbf{w}_{n}}\right) >$ $C(%
\mathbf{u)>0.}$ Now, using the splitting property and (\ref{figa}), we have
that%
\begin{eqnarray*}
\Lambda \left( \mathbf{u+\mathbf{w}_{n}}\right) &=&\frac{E\left( \mathbf{u)+}%
E(\mathbf{\mathbf{w}_{n}}\right) +o(1)}{C\left( \mathbf{u)+}C(\mathbf{%
\mathbf{w}_{n}}\right) +o(1)} \\
&\geq &\min \left( \frac{E\left( \mathbf{u}\right) +o(1)}{C\left( \mathbf{u}%
\right) +o(1)},\frac{E\left( \mathbf{\mathbf{w}_{n}}\right) }{C\left( 
\mathbf{\mathbf{w}_{n}}\right) }\right) .
\end{eqnarray*}%
Then, up to a subsequence , we get (\ref{tiga}).

$\square $

\bigskip

Now we set%
\begin{equation}
\Phi (\mathbf{u})=E(\mathbf{u})+2aC(\mathbf{u})^{s}  \label{phi}
\end{equation}

\begin{equation*}
J_{\delta }(\mathbf{u})=\Lambda \left( \mathbf{u}\right) +\delta \Phi (%
\mathbf{u}),\text{ }\delta >0
\end{equation*}%
and define 
\begin{equation}
\delta _{\infty }=\sup \left\{ \delta >0\ |\ \exists \mathbf{v}:\Lambda
\left( \mathbf{v}\right) +\delta \Phi (\mathbf{v})<\Lambda _{0}\ \right\} .
\label{infinit}
\end{equation}%
By (\ref{hh}), we have that $\delta _{\infty }>0.$

\begin{lemma}
\label{scemo}For any $\delta \geq 0,$ $J_{\delta }(\mathbf{u})\geq \frac{%
\delta }{2}\Phi (\mathbf{u})-M_{\delta }$ where%
\begin{equation*}
M_{\delta }=\ -a\underset{t\geq 0}{\ \min }\left( \frac{\delta }{2}%
t^{s}-t^{s-1}\right) .
\end{equation*}
\end{lemma}

\textbf{Proof: }By assumption (EC-3)(i) we get 
\begin{align*}
J_{\delta }(\mathbf{u})& =\frac{E(\mathbf{u})}{C(\mathbf{u})}+\delta \Phi (%
\mathbf{u})\geq -aC(\mathbf{u})^{s-1}+\frac{\delta }{2}\left[ E(\mathbf{u}%
)+2aC(\mathbf{u})^{s}\right] +\frac{\delta }{2}\Phi (\mathbf{u}) \\
& \geq -aC(\mathbf{u})^{s-1}+\frac{\delta }{2}\left[ -aC(\mathbf{u})^{s}+2aC(%
\mathbf{u})^{s}\right] +\frac{\delta }{2}\Phi (\mathbf{u}) \\
& =-aC(\mathbf{u})^{s-1}+\frac{a\delta }{2}C(\mathbf{u})^{s}+\frac{\delta }{2%
}\Phi (\mathbf{u})\geq \frac{\delta }{2}\Phi (\mathbf{u})-M_{\delta }
\end{align*}%
where%
\begin{equation*}
M_{\delta }=\ -a\underset{t\geq 0}{\ \min }\left( \frac{\delta }{2}%
t^{s}-t^{s-1}\right) .
\end{equation*}

$\square $

\bigskip

\begin{theorem}
\label{legicomp2}Assume that $E\ $and $C$ satisfy (EC-0),...,(EC-3) and the
hylomorphy condition (\ref{hh}). Then, for every $\delta \in \left( 0,\delta
_{\infty }\right) $ (see (\ref{infinit}))$,$ $J_{\delta }$ is $G$-compact
and it has a minimizer $\mathbf{u}_{\delta }\neq 0.$ Moreover $\mathbf{u}%
_{\delta }$ is a minimizer of $E$ on $\mathfrak{M}_{\delta }:=\left\{ 
\mathbf{u}\in X\ |\ C(\mathbf{u})=c_{\delta }\right\} $ where $c_{\delta }=C(%
\mathbf{u}_{\delta }).$
\end{theorem}

\textbf{Proof.} Let $\delta \in \left( 0,\delta _{\infty }\right) ,$ where $%
\delta _{\infty }$ is defined in (\ref{infinit}), and set%
\begin{equation*}
j_{\delta }:=\ \underset{\mathbf{u}\in X}{\inf }J_{\delta }\left( \mathbf{u}%
\right)
\end{equation*}%
By lemma \ref{scemo} and since $\Phi (\mathbf{u})\geq 0$, we have $j_{\delta
}>-\infty $. Then, since $\delta \in \left( 0,\delta _{\infty }\right) ,$ we
have 
\begin{equation}
-\infty <j_{\delta }<\Lambda _{0}  \label{mm}
\end{equation}

Now let $\mathbf{u}_{n}$ be a minimizing sequence of $J_{\delta }$. Let us
prove that $\mathbf{u}_{n}$ is $G$-compact. To this end we shall first prove
that%
\begin{equation*}
\mathbf{u}_{n}\text{ is bounded}.
\end{equation*}%
Arguing by contradiction assume that, up to a subsequence, $\left\Vert 
\mathbf{u}_{n}\right\Vert \longrightarrow +\infty .$ Then, by (EC-3)(ii), we
have 
\begin{equation}
\Phi (\mathbf{u}_{n})=E(\mathbf{u}_{n})+2aC(\mathbf{u}_{n})^{s}%
\longrightarrow +\infty .  \label{ne}
\end{equation}%
By Lemma \ref{scemo} and (\ref{ne}) we get%
\begin{equation*}
J_{\delta }\left( \mathbf{u}_{n}\right) \longrightarrow +\infty .
\end{equation*}%
This contradicts the fact that $\mathbf{u}_{n}$ is a minimizing sequence of $%
J_{\delta }$ and hence $\mathbf{u}_{n}$ is bounded.

Let us prove that 
\begin{equation*}
\mathbf{u}_{n}\text{ is not vanishing.}
\end{equation*}%
By (\ref{mm}) and since $\mathbf{u}_{n}$ is a minimizing sequence for $%
J_{\delta },$ for large $n$ we have 
\begin{equation}
\Lambda \left( \mathbf{u}_{n}\right) \leq J_{\delta }\left( \mathbf{u}%
_{n}\right) <\Lambda _{0}-\eta ,\text{ }\eta >0.  \label{pap}
\end{equation}

Then, by definition of $\Lambda _{0},$ $\mathbf{u}_{n}$ is a not a vanishing
sequence. Hence, by Def. \ref{na}, we can extract a subsequence $\mathbf{u}%
_{n_{k}}$ and we can take a sequence $g_{k}\subset G$ such that $\mathbf{u}%
_{k}^{\prime }:=g_{k}\mathbf{u}_{n_{k}}$ is weakly convergent to some 
\begin{equation}
\mathbf{u}_{\delta }\neq 0.  \label{fi2}
\end{equation}%
We can write 
\begin{equation*}
\mathbf{u}_{n}^{\prime }=\mathbf{u}_{\delta }+\mathbf{w}_{n}
\end{equation*}%
with $\mathbf{w}_{n}\rightharpoonup 0$ weakly. In order to show that $%
J_{\delta }$ is $G$-compact we need to prove that, up to a subsequence, we
have 
\begin{equation*}
\mathbf{w}_{n}\rightarrow 0\text{ strongly}
\end{equation*}%
Clearly we can assume that $\mathbf{w}_{n}\neq 0$ for all $n.$

By the splitting property of $E$ and $C$ and lemma (\ref{leman}), we have
that%
\begin{eqnarray}
j_{\delta } &=&\lim J_{\delta }\left( \mathbf{u}_{\delta }+\mathbf{w}%
_{n}\right) =\lim \left[ \Lambda \left( \mathbf{u}_{\delta }+\mathbf{w}%
_{n}\right) +\delta \Phi \left( \mathbf{u}_{\delta }+\mathbf{w}_{n}\right) %
\right] \geq \\
&\geq &\min \left\{ \Lambda \left( \mathbf{u}_{\delta }\right) ,\lim \Lambda
\left( \mathbf{w}_{n}\right) \right\} +\delta \lim \Phi (\mathbf{u}_{\delta
}+\mathbf{w}_{n})  \label{qu}
\end{eqnarray}

By the splitting property (EC-2) and since $s\geq 1,$ we have that 
\begin{align}
\lim \Phi (\mathbf{u}_{\delta }+\mathbf{w}_{n})& =\lim \left( E(\mathbf{u}%
_{\delta }+\mathbf{w}_{n})+aC(\mathbf{u}_{\delta }+\mathbf{w}_{n})^{s}\right)
\notag \\
& =E(\mathbf{u}_{\delta })+\lim E(\mathbf{w}_{n})+a\lim \left( C(\mathbf{u}%
_{\delta })+C(\mathbf{w}_{n})\right) ^{s}  \notag \\
& \geq E(\mathbf{u}_{\delta })+\lim E(\mathbf{w}_{n})+a\lim \left( C(\mathbf{%
u}_{\delta })^{s}+C(\mathbf{w}_{n})^{s}\right)  \notag \\
& =E(\mathbf{u}_{\delta })+aC(\mathbf{u}_{\delta })^{s}+\lim E(\mathbf{w}%
_{n})+a\lim C(\mathbf{w}_{n})^{s}  \notag \\
& =\Phi (\mathbf{u}_{\delta })+\lim \Phi (\mathbf{w}_{n}).  \label{qw}
\end{align}%
Then by (\ref{qu}) and by (\ref{qw}) we have 
\begin{equation}
j_{\delta }\geq \min \left\{ \Lambda \left( \mathbf{u}_{\delta }\right)
,\lim \Lambda \left( \mathbf{w}_{n}\right) \right\} +\delta \Phi (\mathbf{u}%
_{\delta })+\delta \lim \Phi (\mathbf{w}_{n}).  \label{per}
\end{equation}

Now there are two possibilities: 
\begin{eqnarray*}
\text{(a) }\min \left\{ \Lambda \left( \mathbf{u}_{\delta }\right) ,\lim
\Lambda \left( \mathbf{w}_{n}\right) \right\} &=&\lim \Lambda \left( \mathbf{%
w}_{n}\right) , \\
\text{(b)}\ \min \left\{ \Lambda \left( \mathbf{u}_{\delta }\right) ,\lim
\Lambda \left( \mathbf{w}_{n}\right) \right\} &=&\Lambda \left( \mathbf{u}%
_{\delta }\right) .
\end{eqnarray*}%
\ We will show that the possibility (a) cannot occur. In fact, if it holds,
by (\ref{per}), we have that%
\begin{eqnarray*}
j_{\delta } &\geq &\lim \Lambda \left( \mathbf{w}_{n}\right) +\delta \Phi
\left( \mathbf{u}_{\delta }\right) +\delta \lim \Phi \left( \mathbf{w}%
_{n}\right) \\
&=&\lim J_{\delta }\left( \mathbf{w}_{n}\right) +\delta \Phi \left( \mathbf{u%
}_{\delta }\right) \\
&\geq &j_{\delta }+\delta \Phi \left( \mathbf{u}_{\delta }\right)
\end{eqnarray*}%
and hence, we get that $\Phi \left( \mathbf{u}_{\delta }\right) \leq 0;$
this, by (EC-3)(i), implies that $\mathbf{u}_{\delta }=0,$ contradicting (%
\ref{fi2}). Then the possibility (b) holds and, by (\ref{per}), we have that%
\begin{eqnarray*}
j_{\delta } &\geq &\Lambda \left( \mathbf{u}_{\delta }\right) +\delta \Phi
\left( \mathbf{u}_{\delta }\right) +\delta \lim \Phi \left( \mathbf{w}%
_{n}\right) \\
&=&J_{\delta }\left( \mathbf{u}_{\delta }\right) +\delta \lim \Phi \left( 
\mathbf{w}_{n}\right) \\
&\geq &j_{\delta }+\delta \lim \Phi \left( \mathbf{w}_{n}\right) .
\end{eqnarray*}%
Then, $\lim \Phi \left( \mathbf{w}_{n}\right) \rightarrow 0$ and by
(EC-3)(iii), $\mathbf{w}_{n}\rightarrow 0$ strongly. We conclude that $%
J_{\delta }$ is $G$-compact and $\mathbf{u}_{\delta }$ is a minimizer of $%
J_{\delta }$. Then $\ \mathbf{u}_{\delta }$ minimizes also the functional%
\begin{equation*}
\frac{E(\mathbf{u})}{c_{\delta }}+\delta \left[ E(\mathbf{u})+ac_{\delta
}^{s}\right] =\left( \frac{1}{c_{\delta }}+\delta \right) E(\mathbf{u}%
)+\delta ac_{\delta }^{s}
\end{equation*}%
on the set $\mathfrak{M}_{\delta }=\left\{ \mathbf{u}\in X\ |\ C(\mathbf{u}%
)=c_{\delta }\right\} $ and hence $\mathbf{u}_{\delta }$ minimizes also $E|_{%
\mathfrak{M}_{\delta }}$.

$\square $

\bigskip

In the following $\mathbf{u}_{\delta }$ will denote a minimizer of $%
J_{\delta }.$

\begin{lemma}
\label{MM3}Let the assumptions of Theorem \ref{legicomp2} be satisfied. Let $%
\delta _{1},\delta _{2}\in (0,\delta _{\infty })$ $\delta _{1}<\delta _{2}$
and let $\mathbf{u}_{\delta _{1}},\mathbf{u}_{\delta _{2}}$ be minimizers of 
$J_{\delta _{1}},$ $J_{\delta _{2}}$ respectively. Then the following
inequalities hold:

\begin{itemize}
\item (a) $J_{\delta _{1}}(\mathbf{u}_{\delta _{1}})<J_{\delta _{2}}(\mathbf{%
u}_{\delta _{2}})$

\item (b) $\Phi (\mathbf{u}_{\delta _{1}})\geq \Phi (\mathbf{u}_{\delta
_{2}}),\ $

\item (c) $\Lambda (\mathbf{u}_{\delta _{1}})\leq \Lambda (\mathbf{u}%
_{\delta _{2}}),$

\item (d) $C(\mathbf{u}_{\delta _{1}})\geq C(\mathbf{u}_{\delta _{2}}).$
\end{itemize}
\end{lemma}

\textbf{Proof.} We prove first the inequality ( a) 
\begin{eqnarray*}
J_{\delta _{1}}\left( \mathbf{u}_{\delta _{1}}\right) &=&\Lambda \left( 
\mathbf{u}_{\delta _{1}}\right) +\delta _{1}\Phi (\mathbf{u}_{\delta _{1}})
\\
&\leq &\Lambda \left( \mathbf{u}_{\delta _{2}}\right) +\delta _{1}\Phi (%
\mathbf{u}_{\delta _{2}})\ \ (\text{since }\mathbf{u}_{\delta _{1}}\ \text{%
minimizes }J_{\delta }\text{)} \\
&<&\Lambda \left( \mathbf{u}_{\delta _{2}}\right) +\delta _{2}\Phi (\mathbf{u%
}_{\delta _{2}})\ \ (\text{since }\Phi \ \text{is positive and }\delta
_{1}<\delta _{2}\text{)} \\
&=&J_{\delta _{2}}\left( \mathbf{u}_{\delta _{2}}\right) .
\end{eqnarray*}%
In order to prove inequalities (b) and (c) we set 
\begin{eqnarray*}
\Lambda (\mathbf{u}_{\delta _{1}}) &=&\Lambda (\mathbf{u}_{\delta _{2}})+a \\
\Phi (\mathbf{u}_{\delta _{1}}) &=&\Phi (\mathbf{u}_{\delta _{2}})+b
\end{eqnarray*}%
We need to prove that $b\geq 0\ $and $a\leq 0.$ We have%
\begin{eqnarray}
J_{\delta _{2}}\left( \mathbf{u}_{\delta _{2}}\right) &\leq &J_{\delta _{2}}(%
\mathbf{u}_{\delta _{1}})\Rightarrow  \notag \\
\Lambda \left( \mathbf{u}_{\delta _{2}}\right) +\delta _{2}\Phi (\mathbf{u}%
_{\delta _{2}}) &\leq &\Lambda \left( \mathbf{u}_{\delta _{1}}\right)
+\delta _{2}\Phi (\mathbf{u}_{\delta _{1}})\Rightarrow  \notag \\
\Lambda \left( \mathbf{u}_{\delta _{2}}\right) +\delta _{2}\Phi (\mathbf{u}%
_{\delta _{2}}) &\leq &\left( \Lambda \left( \mathbf{u}_{\delta _{2}}\right)
+a\right) +\delta _{2}\left( \Phi (\mathbf{u}_{\delta _{2}})+b\right)
\Rightarrow  \notag \\
0 &\leq &a+\delta _{2}b.  \label{blu3}
\end{eqnarray}

On the other hand,%
\begin{eqnarray}
J_{\delta _{1}}\left( \mathbf{u}_{\delta _{2}}\right) &\geq &J_{\delta _{1}}(%
\mathbf{u}_{\delta _{1}})\Rightarrow  \notag \\
\Lambda \left( \mathbf{u}_{\delta _{2}}\right) +\delta _{1}\Phi (\mathbf{u}%
_{\delta _{2}}) &\geq &\Lambda \left( \mathbf{u}_{\delta _{1}}\right)
+\delta _{1}\Phi (\mathbf{u}_{\delta _{1}})\Rightarrow  \notag \\
\Lambda \left( \mathbf{u}_{\delta _{2}}\right) +\delta _{1}\Phi (\mathbf{u}%
_{\delta _{2}}) &\geq &\left( \Lambda \left( \mathbf{u}_{\delta _{2}}\right)
+a\right) +\delta _{1}\left( \Phi (\mathbf{u}_{\delta _{2}})+b\right)
\Rightarrow  \notag \\
0 &\geq &a+\delta _{1}b  \label{verde3}
\end{eqnarray}

From (\ref{blu3}) and (\ref{verde3}) we get%
\begin{equation*}
\left( \delta _{2}-\delta _{1}\right) b\geq 0
\end{equation*}%
and hence $b\geq 0.$

Moreover (\ref{blu3}) and (\ref{verde3}) give also 
\begin{equation*}
\left( \frac{1}{\delta _{2}}-\frac{1}{\delta _{1}}\right) a\geq 0
\end{equation*}%
and hence $a\leq 0.$

Finally we prove inequality (d). Arguing by contradiction we assume that%
\begin{equation}
C(\mathbf{u}_{\delta _{1}})<C(\mathbf{u}_{\delta _{2}}).  \label{gil}
\end{equation}%
Then%
\begin{equation}
aC(\mathbf{u}_{\delta _{1}})^{s}<aC(\mathbf{u}_{\delta _{2}})^{s}.
\label{sum}
\end{equation}%
By (c) and (\ref{gil}) we get 
\begin{equation}
C(\mathbf{u}_{\delta _{1}})\Lambda (\mathbf{u}_{\delta _{1}})<C(\mathbf{u}%
_{\delta _{2}})\Lambda (\mathbf{u}_{\delta _{2}}).  \label{sumo}
\end{equation}%
Taking the sum in (\ref{sum}) and (\ref{sumo}) we get%
\begin{equation*}
\Phi (\mathbf{u}_{\delta _{1}})<\Phi (\mathbf{u}_{\delta _{2}})
\end{equation*}%
and this contradicts (b).

$\square $

\begin{lemma}
\label{BB3} Let the assumptions of Theorem \ref{legicomp2} be satisfied and
assume that also (\ref{poo}) is satisfied. Let $\delta _{1},\delta _{2}\in
(0,\delta _{\infty })$ $\delta _{1}<\delta _{2}$ and let $\mathbf{u}_{\delta
_{1}},\mathbf{u}_{\delta _{2}}$ be non zero minimizers of $J_{\delta _{1}},$ 
$J_{\delta _{2}}$ respectively. The following inequalities hold:

\begin{itemize}
\item (a) $\Phi (\mathbf{u}_{\delta _{1}})>\Phi (\mathbf{u}_{\delta _{2}}),\ 
$

\item (b) $\Lambda (\mathbf{u}_{\delta _{1}})<\Lambda (\mathbf{u}_{\delta
_{2}})\ $

\item (c) $C(\mathbf{u}_{\delta _{1}})>C(\mathbf{u}_{\delta _{2}}).$
\end{itemize}
\end{lemma}

\textbf{Proof:} Let $\delta _{1},\delta _{2}\in \left( 0,\delta _{\infty
}\right) $ $\delta _{1}<\delta _{2}.$ By Lemma \ref{MM3} there exist $%
\mathbf{u}_{\delta _{1}},\mathbf{u}_{\delta _{2}}$ non zero minimizers of $%
J_{\delta _{1}},$ $J_{\delta _{2}}.$

By Lemma \ref{MM3}) we know that $\Phi (\mathbf{u}_{\delta _{1}})>\Phi (%
\mathbf{u}_{\delta _{2}}),$ so in order to prove (a) we need only to show
that $\Phi (\mathbf{u}_{\delta _{1}})\neq \Phi (\mathbf{u}_{\delta _{2}}).$
We argue indirectly and assume that 
\begin{equation}
\Phi (\mathbf{u}_{\delta _{1}})=\Phi (\mathbf{u}_{\delta _{2}}).
\label{bla3}
\end{equation}%
By the previous lemma, we have that%
\begin{equation}
\Lambda \left( \mathbf{u}_{\delta _{1}}\right) \leq \Lambda \left( \mathbf{u}%
_{\delta _{2}}\right)  \label{ble3}
\end{equation}%
Also, we have that%
\begin{eqnarray*}
\Lambda \left( \mathbf{u}_{\delta _{2}}\right) +\delta _{2}\Phi \left( 
\mathbf{u}_{\delta _{2}}\right) &\leq &\Lambda \left( \mathbf{u}_{\delta
_{1}}\right) +\delta _{2}\Phi (\mathbf{u}_{\delta _{1}})\ \ \text{(since\ }%
\mathbf{u}_{\delta _{2}}\ \text{minimizes }J_{\delta _{2}}\text{)} \\
&=&\Lambda \left( \mathbf{u}_{\delta _{1}}\right) +\delta _{2}\Phi \left( 
\mathbf{u}_{\delta _{2}}\right) \ \ \text{(by\ (\ref{bla3}))}
\end{eqnarray*}%
and so%
\begin{equation*}
\Lambda \left( \mathbf{u}_{\delta _{2}}\right) \leq \Lambda \left( \mathbf{u}%
_{\delta _{1}}\right)
\end{equation*}%
and by (\ref{ble3}) we get 
\begin{equation}
\Lambda \left( \mathbf{u}_{\delta _{1}}\right) =\Lambda \left( \mathbf{u}%
_{\delta _{2}}\right) .  \label{bli3}
\end{equation}

Then, it follows that $\mathbf{u}_{\delta _{1}}$ is also a minimizer of $%
J_{\delta _{2}};\ $in fact, by (\ref{bli3}) and (\ref{bla3}))%
\begin{eqnarray*}
J_{\delta _{2}}\left( \mathbf{u}_{\delta _{1}}\right) &=&\Lambda \left( 
\mathbf{u}_{\delta _{1}}\right) +\delta _{2}\Phi \left( \mathbf{u}_{\delta
_{1}}\right) \\
&=&\Lambda \left( \mathbf{u}_{\delta _{2}}\right) +\delta _{2}\Phi \left( 
\mathbf{u}_{\delta _{2}}\right) =J_{\delta _{2}}\left( \mathbf{u}_{\delta
_{2}}\right) .
\end{eqnarray*}%
Then, we have that $J_{\delta _{2}}^{\prime }\left( \mathbf{u}_{\delta
_{1}}\right) =0$ as well as $J_{\delta _{1}}\left( \mathbf{u}_{\delta
_{1}}\right) =0$ which esplicitely give%
\begin{eqnarray*}
\Lambda ^{\prime }\left( \mathbf{u}_{\delta _{1}}\right) +\delta _{2}\Phi
^{\prime }\left( \mathbf{u}_{\delta _{1}}\right) &=&0 \\
\Lambda ^{\prime }\left( \mathbf{u}_{\delta _{1}}\right) +\delta _{1}\Phi
^{\prime }\left( \mathbf{u}_{\delta _{1}}\right) &=&0.
\end{eqnarray*}%
The above equations imply that 
\begin{eqnarray*}
\Phi ^{\prime }\left( \mathbf{u}_{\delta _{1}}\right) &=&0 \\
\Lambda ^{\prime }\left( \mathbf{u}_{\delta _{1}}\right) &=&0.
\end{eqnarray*}%
Since $\Lambda \left( \mathbf{u}\right) =\frac{E\left( \mathbf{u}\right) }{%
C\left( \mathbf{u}\right) }$ and $\Phi (\mathbf{u})=E(\mathbf{u})+2aC(%
\mathbf{u})^{s},$ the above system of equations becomes%
\begin{eqnarray}
\frac{E^{\prime }\left( \mathbf{u}_{\delta _{1}}\right) }{C\left( \mathbf{u}%
_{\delta _{1}}\right) }-\frac{E\left( \mathbf{u}_{\delta _{1}}\right) }{%
C\left( \mathbf{u}_{\delta _{1}}\right) ^{2}}C^{\prime }\left( \mathbf{u}%
_{\delta _{1}}\right) &=&0  \notag \\
E^{\prime }(\mathbf{u}_{\delta _{1}})+2asC(\mathbf{u}_{\delta
_{1}})^{s-1}C^{\prime }\left( \mathbf{u}_{\delta _{1}}\right) &=&0.
\label{sega}
\end{eqnarray}%
Eliminating $E^{\prime }(\mathbf{u}_{\delta _{1}}),$ we get%
\begin{equation*}
\left( 2asC(\mathbf{u}_{\delta _{1}})^{s}+E(\mathbf{u}_{\delta _{1}})\right) 
\frac{C^{\prime }\left( \mathbf{u}_{\delta _{1}}\right) }{C\left( \mathbf{u}%
_{\delta _{1}}\right) ^{2}}=0
\end{equation*}%
and, using (\ref{phi}), we get%
\begin{equation*}
\frac{\Phi \left( \mathbf{u}_{\delta _{1}}\right) +2a(s-1)C(\mathbf{u}%
_{\delta _{1}})^{s}}{C\left( \mathbf{u}_{\delta _{1}}\right) ^{2}}C^{\prime
}\left( \mathbf{u}_{\delta _{1}}\right) =0.
\end{equation*}%
By assumption (EC-3) (i) and since $s>1$, we have 
\begin{equation*}
\frac{\Phi \left( \mathbf{u}_{\delta _{1}}\right) +2a(s-1)C(\mathbf{u}%
_{\delta _{1}})^{s}}{C\left( \mathbf{u}_{\delta _{1}}\right) ^{2}}>0,
\end{equation*}%
then $C^{\prime }\left( \mathbf{u}_{\delta _{1}}\right) =0$, and hence, by (%
\ref{sega}), also $E^{\prime }\left( \mathbf{u}_{\delta _{1}}\right) =0.$
Finally by (\ref{poo}) $\mathbf{u}_{\delta _{1}}=0,$ and we get a
contradiction.

In order to prove (b) we argue indirectly and assume that 
\begin{equation}
\Lambda (\mathbf{u}_{\delta _{1}})=\Lambda (\mathbf{u}_{\delta _{2}}).
\label{bla4}
\end{equation}

By (a), we have that%
\begin{equation}
\Phi \left( \mathbf{u}_{\delta _{1}}\right) >\Phi \left( \mathbf{u}_{\delta
_{2}}\right) .  \label{ble4}
\end{equation}

Also, we have that%
\begin{eqnarray*}
\Lambda \left( \mathbf{u}_{\delta _{1}}\right) +\delta _{1}\Phi \left( 
\mathbf{u}_{\delta _{1}}\right) &\leq &\Lambda \left( \mathbf{u}_{\delta
_{2}}\right) +\delta _{1}\Phi (\mathbf{u}_{\delta _{2}})\ \ \text{(since\ }%
\mathbf{u}_{\delta _{1}}\ \text{minimizes }J_{\delta _{1}}\text{)} \\
&\leq &\Lambda \left( \mathbf{u}_{\delta _{1}}\right) +\delta _{1}\Phi
\left( \mathbf{u}_{\delta _{2}}\right) \ \ \text{(by\ (\ref{bla4}))}
\end{eqnarray*}%
and so%
\begin{equation*}
\Phi \left( \mathbf{u}_{\delta _{1}}\right) \leq \Phi \left( \mathbf{u}%
_{\delta _{2}}\right)
\end{equation*}%
and this contadicts (\ref{ble4}).

Let us prove the inequality (c).

Since%
\begin{equation*}
\Lambda \left( \mathbf{u}_{\delta _{i}}\right) C\left( \mathbf{u}_{\delta
_{i}}\right) =E\left( \mathbf{u}_{\delta _{i}}\right) ,\text{ }i=1,2
\end{equation*}%
we have 
\begin{equation*}
\Phi \left( \mathbf{u}_{\delta _{i}}\right) =\Lambda \left( \mathbf{u}%
_{_{\delta _{i}}}\right) C\left( \mathbf{u}_{\delta _{i}}\right) +2aC\left( 
\mathbf{u}_{\delta _{i}}\right) ^{s}\text{, }i=1,2
\end{equation*}%
and the conclusion easily follows from inequalities (a) and (b).

$\square $\bigskip

\subsection{The stability result\label{35}}

In the previous subsection \ref{33}, we have proved the existence of
minimizers, namely that $\Gamma (e,c)\neq \varnothing $ (see (\ref{plis})).
In this subsection, we prove the stability of $\Gamma (e,c)$ namely that the
minimizers are hylomorphic solitons. More exactly we will prove the
following two theorems:

\begin{theorem}
\label{astra1}Assume that $E\ $and $C$ satisfy (EC-0),...,(EC-2),\textit{\
(EC-3). Assume also that }the hylomorphy condition of Def. \ref{dhc} is
satisfied. Then for any $\delta \in \left( 0,\delta _{\infty }\right) $ ($%
\delta _{\infty }>0$ defined in (\ref{infinit})) there exists a hylomorphic
soliton $\mathbf{u}_{\delta }.$ \textit{Moreover assume that} 
\begin{equation}
\left\Vert E^{\prime }(\mathbf{u})\right\Vert +\left\Vert C^{\prime }(%
\mathbf{u})\right\Vert =0\Leftrightarrow \mathbf{u}=0.  \label{poo}
\end{equation}

Then, if $\delta _{1}<\delta _{2},$ the corresponding solitons $\mathbf{u}%
_{\delta _{1}},\mathbf{u}_{\delta _{2}}$ are distinct, and we have that

\begin{itemize}
\item (a) $\Lambda (\mathbf{u}_{\delta _{1}})<\Lambda (\mathbf{u}_{\delta
_{2}})$

\item (b) $C(\mathbf{u}_{\delta _{1}})>C(\mathbf{u}_{\delta _{2}}).$

\item (c) $E(\mathbf{u}_{\delta _{1}})+aC(\mathbf{u}_{\delta _{1}})^{s}>E(%
\mathbf{u}_{\delta _{2}})+aC(\mathbf{u}_{\delta _{2}})^{s}$
\end{itemize}
\end{theorem}

\bigskip

\begin{rem}
Variants of the above results have been stated in \cite{befolak} and \cite%
{befocao}.
\end{rem}

Before proving Theorem \ref{astra1} we need to recall some result.

\begin{theorem}
\label{propV}Let $\Gamma $ be an invariant set and assume that there exists
a differentiable real function $V$ (called a Liapunov function) defined on a
neighborhood of $\Gamma $ such that

\begin{itemize}
\item (a) $V(\mathbf{u})\geq 0$ and\ $V(\mathbf{u})=0\Leftrightarrow u\in
\Gamma $

\item (b) $\partial _{t}V(\gamma _{t}\left( \mathbf{u}\right) )\leq 0$

\item (c) $V(\mathbf{u}_{n})\rightarrow 0\Leftrightarrow d(\mathbf{u}%
_{n},\Gamma )\rightarrow 0.$
\end{itemize}

\noindent Then $\Gamma$ is stable.
\end{theorem}

\textbf{Proof. }This is a classical result. A proof of it in this form can
be found in \cite{befolak} or \cite{befocao}.

$\square $

We shall need also the following Lemma

\begin{lemma}
\label{LC}Let $V\geq 0$ be $G$-compact functional and let $\Gamma =V^{-1}(0)$
be the set of minimizers of $V.\ $If $\Gamma \neq \varnothing ,$ then $%
\Gamma $ is $G$-compact and $V$ satisfies the point (c) of the previous
lemma.
\end{lemma}

\textbf{Proof}: A proof can be found in \cite{befolak} or \cite{befocao}.

$\square$

\textbf{Proof of Th. \ref{astra1}.} Let $\mathbf{u}_{\delta }$ be a
minimizer of $E$ on

\begin{equation*}
\mathfrak{M}_{\delta }=\left\{ \mathbf{u}\in X\ |\ C(\mathbf{u})=c_{\delta
}\right\}
\end{equation*}
as in Theorem \ref{legicomp2}. It remains to show that 
\begin{equation*}
\Gamma \left( e_{\delta },c_{\delta }\right) =\left\{ \mathbf{u}\in X\ |\ C(%
\mathbf{u})=c_{\delta },E(\mathbf{u})=e_{\delta }\right\} ,\text{ }%
(e_{\delta }=E(\mathbf{u}_{\delta }))
\end{equation*}%
is $G$-compact and stable.

\begin{itemize}
\item $\Gamma \left( e_{\delta },c_{\delta }\right) $ is $G$-compact.
\end{itemize}

To this end, by Lemma \ref{LC}, it will be enough to show that

\begin{equation*}
V(\mathbf{u})=(E(\mathbf{u})-e_{\delta })^{2}+(C(\mathbf{u})-c_{\delta })^{2}
\end{equation*}%
is $G$ compact.

Let $\mathbf{w}_{n}$ be a minimizing sequence for $V,$ then $V\left( \mathbf{%
w}_{n}\right) \rightarrow 0$ and consequently $E\left( \mathbf{w}_{n}\right)
\rightarrow e_{\delta }$ and $C\left( \mathbf{w}_{n}\right) \rightarrow
c_{\delta }$. Now, since 
\begin{equation*}
\inf J_{\delta }=\frac{e_{\delta }}{c_{\delta }}+\delta \left[ e_{\delta
}+ac_{\delta }^{s}\right] ,
\end{equation*}%
we have that $\mathbf{w}_{n}$ is a minimizing sequence also for $J_{\delta
}. $ Then, since by Theorem \ref{legicomp2} $J_{\delta }$ is $G$-compact, we
get that 
\begin{equation}
\mathbf{w}_{n}\ \text{is}\ G\text{-compact}.  \label{paracula}
\end{equation}%
So we conclude that $V$ is $G$-compact.

\begin{itemize}
\item $\Gamma \left( e_{\delta },c_{\delta }\right) $ is stable.
\end{itemize}

In fact, since $V$ is $G$-compact, by Lemma \ref{LC} we deduce that $%
V^{-1}(0)$ $=\Gamma \left( e_{\delta },c_{\delta }\right) $ satisfies the
point (c) in Theorem \ref{propV}. Moreover clearly $V$ satisfies also the
points (a) and (b) in Theorem \ref{propV}. So, by Theorem \ref{propV}, we
conclude that $\Gamma \left( e_{\delta },c_{\delta }\right) $ is stable.

Finally, if we assume (\ref{poo}), we can use Lemma \ref{BB3} to get
different solitons for different values of $\delta .$ Namely for $\delta
_{1}<\delta _{2}$ we have $\Lambda (\mathbf{u}_{\delta _{1}})<\Lambda (%
\mathbf{u}_{\delta _{2}})$ and $C(\mathbf{u}_{\delta _{1}})>C(\mathbf{u}%
_{\delta _{2}})$.

$\square $

\section{The nonlinear Schr\H{o}dinger Maxwell equation\label{SSE}}

In this section we derive a system of equations (NSM) obtained by coupling
the nonlinear Schr\"{o}dinger equation with Maxwell equations and we prove
the existence of a family of stable solitary waves.

\subsection{General features\label{ec}}

The Schr\H{o}dinger equation for a particle which moves in a potential $V(x)$
is given by 
\begin{equation*}
i\frac{\partial \psi }{\partial t}=-\frac{1}{2}\Delta \psi +V(x)\psi
\end{equation*}%
where $\psi :\mathbb{R\times R}^{3}\rightarrow \mathbb{C\ }$and $V:\mathbb{R}%
^{3}\rightarrow \mathbb{R}$.

We are interested to the nonlinear Schr\H{o}dinger equation: 
\begin{equation}
i\frac{\partial \psi }{\partial t}=-\frac{1}{2}\Delta \psi +\frac{1}{2}%
W^{\prime }(\psi )+V(x)\psi  \label{NSV}
\end{equation}%
where $W:\mathbb{C\rightarrow R}$ and 
\begin{equation}
W^{\prime }(\psi )=\frac{\partial W}{\partial \psi _{1}}+i\frac{\partial W}{%
\partial \psi _{2}}.  \label{w'}
\end{equation}%
We assume that $W$ depends only on $\left\vert \psi \right\vert $, namely 
\begin{equation*}
W(\psi )=F(\left\vert \psi \right\vert )\ \text{and so\ }W^{\prime }(\psi
)=F^{\prime }(\left\vert \psi \right\vert )\frac{\psi }{\left\vert \psi
\right\vert }.
\end{equation*}%
for some smooth function $F:\left[ 0,\infty \right) \rightarrow \mathbb{R}.$
In the following we shall identify, with some abuse of notation, $W$ with $%
F. $

If $V(x)=0,$ then we get the equation 
\begin{equation}
i\frac{\partial \psi }{\partial t}=-\frac{1}{2}\Delta \psi +\frac{1}{2}%
W^{\prime }(\psi );  \tag{NS}  \label{NS}
\end{equation}

Equation (\ref{NSV}) is the Euler-Lagrange equation relative to the
Lagrangian density 
\begin{equation}
\mathcal{L}_{s}=\func{Re}\left( i\partial _{t}\psi \overline{\psi }\right) -%
\frac{1}{2}\left\vert \nabla \psi \right\vert ^{2}-W\left( \psi \right)
-V(x)\left\vert \psi \right\vert ^{2}  \label{quar}
\end{equation}

Now we want to couple the Schr\"{o}dinger equation with the Maxwell
equations. We recall that the use of the covariant derivative provides a
very elegant procedure to combine relativistic field equations (Dirac,
Klein-Gordon etc.) with the Maxwell equations (see e.g. \cite{yangL}, \cite%
{rub}, \cite{befogranas}, \cite{befo11max}). It is possible to use this
procedure also to couple Schr\H{o}dinger and Maxwell equations. This
situation describes the interaction between a charged "matter field" with
the electromagnetic field when the relativistic effects are negligible (see 
\cite{sanchez} and its references).

Let us see how this procedure works. We denote by $\mathbf{E}$, $\mathbf{H}$
the electric and the magnetic field and by $\varphi :\mathbb{R}%
^{3}\rightarrow \mathbb{R}$ and $\mathbf{A}:\mathbb{R}^{3}\rightarrow 
\mathbb{R}^{3},$ $\mathbf{A}=(A_{1},A_{2},A_{3})$ their gauge potentials,
namely fields such that 
\begin{equation*}
\mathbf{E=-}\frac{\partial \mathbf{A}}{\partial t}-\nabla \varphi ,\text{ }%
\mathbf{H=\nabla \times A.}
\end{equation*}%
Now couple (\ref{NSV}) with Maxwell equations by means of the covariant
derivatives. So $\mathcal{L}_{s}$ becomes%
\begin{equation*}
\mathcal{L}_{c}=\func{Re}\left( iD_{t}\psi \overline{\psi }\right) -\frac{1}{%
2}\left\vert \mathbf{D}_{x}\psi \right\vert ^{2}-W\left( \psi \right)
-V(x)\left\vert \psi \right\vert ^{2},
\end{equation*}%
where $\mathbf{D}_{x},\ D_{t}$ denote the covariant derivatives 
\begin{equation*}
\mathbf{D}_{x}\psi =\left( D_{1}\psi ,D_{2}\psi ,D_{3}\psi \right)
\end{equation*}
\begin{equation*}
D_{t}=\frac{\partial }{\partial t}+iq\varphi ,\ D_{j}=\frac{\partial }{%
\partial x^{j}}-iqA_{j}
\end{equation*}%
and $q$ denotes a positive coupling constant wich represents the "strenght"
of the interaction. Adding to $\mathcal{L}_{c}$ the Lagrangian related to
the Maxwell equations 
\begin{equation*}
\mathcal{L}_{M}=\frac{1}{2}\left( \left\vert \frac{\partial \mathbf{%
\QTR{mathbf}{A}}}{\partial t}+\nabla \varphi \right\vert ^{2}-\left\vert
\nabla \times \mathbf{A}\right\vert ^{2}\right) ,
\end{equation*}%
we get the total Lagrangian 
\begin{equation}
\mathcal{L}=\mathcal{L}_{c}+\mathcal{L}_{M}.  \label{total}
\end{equation}

So the total action is 
\begin{equation}
\mathcal{S}=\int \mathcal{L}dxdt.  \label{vetrina}
\end{equation}

If we write $\psi $ in polar form 
\begin{equation*}
\psi (x,t)=u(x,t)\,e^{iS(x,t)},\;\;u\geq 0,\;\;S\in \mathbb{R}/2\pi \mathbb{Z%
}
\end{equation*}%
the action (\ref{vetrina}) takes the following form 
\begin{align}
\mathcal{S(}u,S,\varphi ,\mathbf{A})& =-\int \int \left[ \frac{1}{2}%
\left\vert \nabla u\right\vert ^{2}+V(x)u^{2}+W(u)\right] dxdt+  \notag \\
& -\int \int \left[ \left( \frac{\partial S}{\partial t}+q\varphi \right) +%
\frac{1}{2}\left\vert \nabla S-q\mathbf{A}\right\vert ^{2}\right] \,u^{2}dxdt
\label{azion} \\
& +\frac{1}{2}\int \int \left( \left\vert \frac{\partial \mathbf{A}}{%
\partial t}\mathbf{+}\nabla \varphi \right\vert ^{2}-\left\vert \nabla
\times \mathbf{A}\right\vert ^{2}\right) dxdt.  \notag
\end{align}

Making the variations of $\mathcal{S}$ with respect $u,S,\varphi ,\mathbf{A}$
we get respectively the equations : 
\begin{equation}
-\frac{1}{2}\Delta u+\frac{1}{2}W^{\prime }(u)+\left[ \frac{1}{2}\left\vert
\nabla S-q\mathbf{A}\right\vert ^{2}+\left( \frac{\partial S}{\partial t}%
+q\varphi +V(x)\right) \right] \,u=0  \label{f1}
\end{equation}%
\begin{equation}
\frac{\partial (u^{2})}{\partial t}+\nabla \cdot \left[ \left( \nabla S-q%
\mathbf{A}\right) u^{2}\right] =0  \label{f2}
\end{equation}

\begin{equation}
-\nabla \cdot \left( \frac{\partial \mathbf{\QTR{mathbf}{A}}}{\partial t}%
+\nabla \varphi \right) =qu^{2}\;  \label{f3}
\end{equation}%
\begin{equation}
\nabla \times \left( \nabla \times \mathbf{A}\right) +\frac{\partial }{%
\partial t}\left( \frac{\partial \mathbf{\QTR{mathbf}{A}}}{\partial t}%
+\nabla \varphi \right) =q\left( \nabla S-q\mathbf{A}\right) u^{2}.\;
\label{f4}
\end{equation}%
The last two equations (\ref{f3}) and (\ref{f4}) are the second couple of
the Maxwell equations%
\index{equations!Maxwell} (Gauss and Ampere laws) with respect to a matter
distribution whose electric charge and current density are respectively $%
\rho $ and $\mathbf{j}$ defined by:

\begin{equation}
\rho =qu^{2}  \label{ppp3}
\end{equation}

\begin{equation}
\mathbf{j}=q\left( \nabla S-q\mathbf{A}\right) u^{2}.  \label{p44}
\end{equation}%
Notice that equation (\ref{f2}) is a continuity equation which gives rice to
the conservation of the hylenic charge $C$

\begin{equation}
C=\dint u^{2}.  \label{rec}
\end{equation}%
and hence also to the conservation of the electic charge $qC=q\dint u^{2}.$

Moreover (\ref{f2}) is easily derived from (\ref{f3}) and (\ref{f4}). In
conclusion our system of equations is reduced to (\ref{f1}), (\ref{f3}), (%
\ref{f4}).

Observe that in the electrostatic case. i.e. when 
\begin{equation*}
\frac{\partial u}{\partial t}=0,\text{ }S=\omega t,\text{ }\omega \text{ real%
},\text{ }\mathbf{A}=0,
\end{equation*}%
the system (\ref{f1}), (\ref{f3}), (\ref{f4}) reduces to the system 
\begin{equation}
-\frac{1}{2}\Delta u+\frac{1}{2}W^{\prime }(u)+\left( q\varphi +\omega
+V(x)\right) \,u=0  \label{aaa}
\end{equation}%
\begin{equation}
-\Delta \varphi =qu^{2}\;.  \label{bbb}
\end{equation}%
System (\ref{aaa}), (\ref{bbb}) is called nonlinear Schr\"{o}dinger-Maxwell
system or nonlinear Schr\"{o}dinger-Poisson system and it will be denoted by
NSM.

Observe that, if we consider $\varphi $ as a scalar field (and not the
time-component of a 4-vector) the Schr\"{o}dinger-Poisson equations are
invariant under the Galileo group if $V$ is constant.

Now we compute the energy $E$ related to the system (\ref{f1},...,\ref{f4}).

\begin{theorem}
If $(u,S,\varphi ,\mathbf{A})$ satisfy the Gauss equations (\ref{f3}), the
energy $E$ related to the system (\ref{f1},...,\ref{f4}) takes the following
form:%
\begin{eqnarray*}
E &=&\int \left[ \frac{1}{2}\left\vert \nabla u\right\vert
^{2}+V(x)u^{2}+W(u)\right] dx \\
&&+\frac{1}{2}\int \left[ \left\vert \nabla S-q\mathbf{A}\right\vert
^{2}\,u^{2}\right] dx+\frac{1}{2}\int \left[ \left\vert \frac{\partial 
\mathbf{A}}{\partial t}\mathbf{+}\nabla \varphi \right\vert ^{2}+\left\vert
\nabla \times \mathbf{A}\right\vert ^{2}\right] dx.
\end{eqnarray*}
\end{theorem}

\textbf{Proof.} The Lagrangian $\mathcal{L}$ related to the system (\ref{f1}%
,...,\ref{f4}) is 
\begin{align}
\mathcal{L}& =-\frac{1}{2}\left\vert \nabla u\right\vert ^{2}-V(x)u^{2}-W(u)-
\label{ll} \\
& -\left( \frac{\partial S}{\partial t}+q\varphi \right) u^{2}-\frac{1}{2}%
\left\vert \nabla S-q\mathbf{A}\right\vert ^{2}\,u^{2}  \notag \\
& +\frac{1}{2}\left( \left\vert \frac{\partial \mathbf{A}}{\partial t}%
\mathbf{+}\nabla \varphi \right\vert ^{2}-\left\vert \nabla \times \mathbf{A}%
\right\vert ^{2}\right) .  \notag
\end{align}%
This Lagrangian does not depend on $\frac{\partial u}{\partial t}$ and $%
\frac{\partial \varphi }{\partial t}.$ Then the related energy is (see \cite%
{Gelfand} chapter 7) 
\begin{equation*}
E=\int \left[ \frac{\partial \mathcal{L}}{\partial \left( \frac{\partial S}{%
\partial t}\right) }\cdot \frac{\partial S}{\partial t}+\frac{\partial 
\mathcal{L}}{\partial \left( \frac{\partial \mathbf{A}}{\partial t}\right) }%
\cdot \frac{\partial \mathbf{A}}{\partial t}-\mathcal{L}\right] dx.
\end{equation*}

So, by a direct calculation, we get 
\begin{equation*}
E=\dint \left( \frac{\partial \mathbf{A}}{\partial t}\mathbf{+}\nabla
\varphi \right) \cdot \frac{\partial \mathbf{A}}{\partial t}+q\varphi u^{2}+%
\frac{1}{2}\left( \left\vert \nabla u\right\vert ^{2}+\left\vert \nabla S-q%
\mathbf{A}\right\vert ^{2}\,u^{2}\right) +
\end{equation*}%
\begin{equation}
V(x)u^{2}+W(u)-\frac{1}{2}\left\vert \frac{\partial \mathbf{A}}{\partial t}%
\mathbf{+}\nabla \varphi \right\vert ^{2}+\frac{1}{2}\left\vert \nabla
\times \mathbf{A}\right\vert ^{2}.  \label{energit}
\end{equation}

By the Gauss equation (\ref{f3}), multiplying by $\varphi $ and integrating,
we get 
\begin{equation}
\int q\varphi u^{2}=\int \nabla \varphi \cdot \left( \frac{\partial \mathbf{A%
}}{\partial t}\mathbf{+}\nabla \varphi \right) .  \label{gaussb}
\end{equation}

The above equality (\ref{gaussb}) easily implies that 
\begin{equation*}
\dint q\varphi u^{2}+\left( \frac{\partial \mathbf{A}}{\partial t}\mathbf{+}%
\nabla \varphi \right) \cdot \frac{\partial \mathbf{A}}{\partial t}-\frac{1}{%
2}\left\vert \frac{\partial \mathbf{A}}{\partial t}\mathbf{+}\nabla \varphi
\right\vert ^{2}=
\end{equation*}%
\begin{equation}
\frac{1}{2}\dint \left\vert \frac{\partial \mathbf{A}}{\partial t}\mathbf{+}%
\nabla \varphi \right\vert ^{2}.  \label{inserto}
\end{equation}

Inserting (\ref{inserto}) into (\ref{energit}) we get the conclusion.

$\square $

.

\subsection{Statement of the results}

We make the following assumptions on $W$ and $V$\label{ppagg}

\begin{itemize}
\item \textbf{(W-i)} $W$ is a $C^{2}$ function s.t. 
\begin{equation}
W(0)=W^{\prime }(0)=0\ \text{and }W^{\prime \prime }(0)=2E_{0}>0;
\label{imis}
\end{equation}

\item \textbf{(W-ii)} if we set 
\begin{equation}
W(s)=E_{0}s^{2}+N(s),  \label{W}
\end{equation}%
then%
\begin{equation}
\exists s_{0}\in \mathbb{R}^{+}\text{ such that }N(s_{0})<-V_{0}s_{0}^{2}
\label{W1}
\end{equation}%
where 
\begin{equation*}
V_{0}=\max V;
\end{equation*}

\item \textbf{(W-iii) }there exist $q,$ $r$ in $(2,6),\ \ $ s. t.%
\begin{equation}
|N^{\prime }(s)|\leq c_{1}s^{r-1}+c_{2}s^{q-1}  \label{Wp}
\end{equation}

\item \textbf{(W-iv)}%
\begin{equation}
N(s)\geq -cs^{p},\text{ }c\geq 0,\ 2<p<2+\frac{4}{3}\text{ for }s\text{ large%
}  \label{W0}
\end{equation}
\end{itemize}

$\bigskip $

$V:\mathbb{R}^{3}\rightarrow \mathbb{R}$ being a potential function
satisfying the assumptions:

\begin{itemize}
\item \textbf{(V-i)} $V$ continuous and 
\begin{equation}
V(x)\geq 0,\text{ }x\in \mathbb{R}^{3}  \label{V0}
\end{equation}

\item \textbf{(V-ii)} $V$ is a lattice potential, namely it satisfies the
periodicity condition:%
\begin{equation}
V(x)=V(x+Az)\text{ for all }x\in \mathbb{R}^{3}\text{ and }z\in \mathbb{Z}%
^{3}  \label{V1'}
\end{equation}%
where $A$ is a $3\times 3$ invertible matrix.
\end{itemize}

If we set%
\begin{equation*}
\text{ }\mathbf{E=-}\frac{\partial \mathbf{\QTR{mathbf}{A}}}{\partial t}%
-\nabla \varphi ,\text{ }\mathbf{H}=\nabla \times \mathbf{A},\ \Theta
=\left( \nabla S-q\mathbf{A}\right) u\text{,}
\end{equation*}%
the energy $E$ takes the form%
\begin{equation}
E=\int \left( \frac{1}{2}\left\vert \nabla u\right\vert ^{2}+V(x)u^{2}+W(u)+%
\frac{1}{2}\left( \Theta ^{2}+\mathbf{E}^{2}+\mathbf{H}^{2}\right) \right)
dx.  \label{ens}
\end{equation}%
Instead of using the variables $\left( u,S,\mathbf{E},\mathbf{H}\right) ,$
we will use the variables $\left( u,\Theta ,\mathbf{E},\mathbf{H}\right) $
so that the generic point in the phase space is given by%
\begin{equation*}
\mathbf{u}=\left( u,\Theta ,\mathbf{E},\mathbf{H}\right)
\end{equation*}%
and the phase space is given by%
\begin{equation}
X=\left\{ \mathbf{u}=\left( u,\Theta ,\mathbf{E},\mathbf{H}\right) \in
H^{1}\left( \mathbb{R}^{3}\right) \times L^{2}\left( \mathbb{R}^{3}\right)
^{9}:\nabla \cdot \mathbf{E}=qu^{2}\right\}  \label{face}
\end{equation}%
where $H^{1}\left( \mathbb{R}^{3}\right) $ is the usual Sobolev space.

We equip $X$ with the norm related to the quadratic part of the energy,
namely:%
\begin{equation}
\left\Vert \mathbf{u}\right\Vert ^{2}=\int \left\vert \nabla u\right\vert
^{2}+2E_{0}u^{2}+\Theta ^{2}+\mathbf{E}^{2}+\mathbf{H}^{2}  \label{le}
\end{equation}%
where $E_{0}$ is defined by (\ref{W}). Then the energy $E$ can be written as
follows:%
\begin{equation}
E=\frac{1}{2}\left\Vert \mathbf{u}\right\Vert ^{2}+\int V(x)u^{2}+\int N(u).
\label{enss}
\end{equation}

We notice that the new variables do not change the expression for the
charge, namely $C$ keeps the form (\ref{rec}). Finally, as usual 
\begin{equation*}
\Lambda =\frac{E\left( \mathbf{u}\right) }{C\left( \mathbf{u}\right) }
\end{equation*}%
will denote the hylenic ratio.

In the following we shall assume that the Cauchy problem for the system (\ref%
{f1}, \ref{f3}, \ref{f4}) is well posed in $X$ and we refer to \cite%
{guonastrauss}, \cite{nakats} and \cite{nakawa} for some results in this
direction.

We shall prove the following existence results of hylomorphic solitons for
NSM.

\begin{theorem}
\label{exnsm}Let $W$ and $V$ satisfy the assumptions (WB-i),...,(WB-iv) and
(V-i),(V-ii). Then, if $q>0$ is sufficiently small, there exists $\delta
_{\infty }>0$ such that the dynamical system described by the system (\ref%
{f1}), (\ref{f3}) (\ref{f4}) has a family $\mathbf{u}_{\delta }=\left(
u_{\delta },\Theta _{\delta },\mathbf{E}_{\delta },\mathbf{H}_{\delta
}\right) $ ($\delta \in \left( 0,\delta _{\infty }\right) )$ of hylomorphic
solitons (Definition \ref{tdc} ). Moreover if $\delta _{1}<\delta _{2}$ we
have that

\begin{itemize}
\item (a) $\Lambda (\mathbf{u}_{\delta _{1}})<\Lambda (\mathbf{u}_{\delta
_{2}})$

\item (b) $\left\Vert u_{\delta _{1}}\right\Vert _{L^{2}}>\ \left\Vert
u_{\delta _{2}}\right\Vert _{L^{2}}$
\end{itemize}
\end{theorem}

\begin{theorem}
\label{nsmbis}The solitons $\mathbf{u}_{\delta }=\left( u_{\delta },\Theta
_{\delta },\mathbf{E}_{\delta },\mathbf{H}_{\delta }\right) $ in Theorem \ref%
{exnsm} are stationary solutions of (\ref{f1}), (\ref{f3}) (\ref{f4})), this
means that $\Theta _{\delta }=\mathbf{H}_{\delta }=0,$ $\mathbf{E}_{\delta
}=-\nabla \varphi _{\delta }$, $u_{\delta },\varphi _{\delta }$ do not
depend on $t$ and they solve, for suitable real numbers $\omega ,$ the
nonlinear Schr\H{o}dinger-Poisson system 
\begin{align}
-\frac{1}{2}\Delta u_{\delta }+V(x)u_{\delta }+\frac{1}{2}W^{\prime
}(u_{\delta })+q\varphi _{\delta }u_{\delta }& =-\omega u_{\delta }
\label{din} \\
-\Delta \varphi _{\delta }& =qu_{\delta }^{2}.  \label{dan}
\end{align}

\begin{remark}
If the coupling constant $q=0$ equations (\ref{f1}), (\ref{f3}) (\ref{f4})
reduce to the Schr\"{o}dinger equation and Theorem \ref{exnsm} becomes in
this case a variant of well known stability results (see \cite{Cazli}, \cite%
{sulesulem} and its references).
\end{remark}
\end{theorem}

The proof of Theorem \ref{exnsm} is based on the abstract Theorem \ref%
{astra1+}. First of all observe that, since $V$ satisfies (\ref{V1'}), the
energy $E$ is invariant under the representation $T_{z}$ of the group $G:=%
\mathbb{Z}^{3}$ 
\begin{equation*}
T_{z}\mathbf{u}(x)=\mathbf{u}(x+Az),\ \ \ z\in \mathbb{Z}^{3}
\end{equation*}%
where $A$ is as in (\ref{V1'}).

\subsection{Proof of the results}

In this section we shall prove that $E$ and $C$ satisfy assumptions (EC-2)
(splitting), (EC-3) (coercivity) and the hylomorphy assumption.

\begin{lemma}
\label{split}Let the assumptions of theorem \ref{exnsm} be satisfied. Then $%
E $ and $C,$ defined by (\ref{enss}) and (\ref{rec}) satisfy the splitting
property (EC-2).
\end{lemma}

\textbf{Proof. } For any $\mathbf{u=}\left( u\mathbf{,\Theta ,\mathbf{E},%
\mathbf{H}}\right) \in X$ the energy $E(\mathbf{u})$ in (\ref{enss}) can be
written

\begin{equation*}
E\left( \mathbf{u}\right) =A(\mathbf{u},\mathbf{u})+K(u\mathbf{)}
\end{equation*}

where%
\begin{equation*}
A(\mathbf{u},\mathbf{u})=\frac{1}{2}\left\Vert \mathbf{u}\right\Vert
^{2}+\int V(x)u^{2}
\end{equation*}

and

\begin{equation}
K(u\mathbf{)=}\int N\left( u\right) dx.  \label{cappio}
\end{equation}

The hylenic charge $C(\mathbf{u})=\int u^{2}$ and $A(\mathbf{u},\mathbf{u})$
are quadratic forms, then, by remark \ref{quadratic}, they satisfy the
splitting property. So, in order to show that also the energy $E\left( 
\mathbf{u}\right) $ satisfies (EC-2), we have only to show that $K(u\mathbf{)%
}$ in (\ref{cappio}) satisfies the splitting property. Let $H^{1}(\mathbb{R}%
^{3})$ denote the usual Sobolev space, then for any measurable $A\subset 
\mathbb{R}^{3}$ and any $u\in H^{1}(\mathbb{R}^{3})$, we set%
\begin{equation*}
K_{A}(u\mathbf{)=}\int_{A}N(u)dx.
\end{equation*}%
Now consider any sequence 
\begin{equation*}
u_{n}=u+w_{n}\in H^{1}(\mathbb{R}^{3})
\end{equation*}%
where $w_{n}$ converges weakly to $0.$

Choose $\varepsilon >0$ and $R=R(\varepsilon )>0$ such that 
\begin{equation}
\left\vert K_{B_{R}^{c}}\left( u\right) \right\vert <\varepsilon  \label{bis}
\end{equation}

where $\ $%
\begin{equation*}
B_{R}^{c}=\mathbb{R}^{3}-B_{R}\text{ and }B_{R}=\left\{ x\in \mathbb{R}%
^{3}:\left\vert x\right\vert <R\right\} .
\end{equation*}%
Since $w_{n}\rightharpoonup 0$ weakly in $H^{1}\left( \mathbb{R}^{3}\right) $%
, by usual compactness arguments, we have that 
\begin{equation}
K_{B_{R}}\left( w_{n}\right) \rightarrow 0\text{ and }K_{B_{R}}\left(
u+w_{n}\right) \rightarrow K_{B_{R}}\left( u\right) .  \label{bibis}
\end{equation}

Then, by (\ref{bis}) and (\ref{bibis}), we have

\begin{align}
& \underset{n\rightarrow \infty }{\lim }\left\vert K\left( u+w_{n}\right)
-K\left( u\right) -K\left( w_{n}\right) \right\vert  \notag \\
& =\ \underset{n\rightarrow \infty }{\lim }\left\vert K_{B_{R}^{c}}\left(
u+w_{n}\right) +K_{B_{R}}\left( u+w_{n}\right) -K_{B_{R}^{c}}\left( u\right)
-K_{B_{R}}\left( u\right) -K_{B_{R}^{c}}\left( w_{n}\right) -K_{B_{R}}\left(
w_{n}\right) \right\vert  \notag \\
& \mathbb{=}\ \underset{n\rightarrow \infty }{\lim }\left\vert
K_{B_{R}^{c}}\left( u+w_{n}\right) -K_{B_{R}^{c}}\left( u\right)
-K_{B_{R}^{c}}\left( w_{n}\right) \right\vert  \notag \\
& \leq \ \underset{n\rightarrow \infty }{\lim }\left\vert
K_{B_{R}^{c}}\left( u+w_{n}\right) -K_{B_{R}^{c}}\left( w_{n}\right)
\right\vert +\varepsilon .  \label{pepe}
\end{align}

Now, by the intermediate value theorem, there exists $\zeta _{n}\in (0,1)$
such that for $z_{n}=$ $\zeta _{n}u+\left( 1-\zeta _{n}\right) w_{n}$, we
have that

\begin{align*}
\left\vert K_{B_{R}^{c}}\left( u+w_{n}\right) -K_{B_{R}^{c}}\left(
w_{n}\right) \right\vert & =\left\vert \left\langle K_{B_{R}^{c}}^{\prime
}\left( z_{n}\right) ,u\right\rangle \right\vert \\
& \leq \int_{B_{R}^{c}}\left\vert N^{\prime }(z_{n})u\right\vert \leq (\text{%
by (\ref{Wp})}) \\
& \leq \int_{B_{R}^{c}}c_{1}\left\vert z_{n}\right\vert ^{r-1}\left\vert
u\right\vert +c_{2}\left\vert z_{n}\right\vert ^{q-1}\left\vert u\right\vert
\\
& \leq c_{1}\left\Vert z_{n}\right\Vert _{L^{r}(B_{R}^{c})}^{r-1}\left\Vert
u\right\Vert _{L^{r}(B_{R}^{c})}+c_{2}\left\Vert z_{n}\right\Vert
_{L^{q}(B_{R}^{c})}^{q-1}\left\Vert u\right\Vert _{_{L^{q}(B_{R}^{c})}} \\
& (\text{if }R\text{ is large enough)} \\
& \leq c_{3}\left( \left\Vert z_{n}\right\Vert
_{L^{r}(B_{R}^{c})}^{r-1}+\left\Vert z_{n}\right\Vert
_{L^{q}(B_{R}^{c})}^{q-1}\right) \varepsilon .
\end{align*}%
So we have 
\begin{equation}
\left\vert K_{B_{R}^{c}}\left( u+w_{n}\right) -K_{B_{R}^{c}}\left(
w_{n}\right) \right\vert \leq c_{3}\left( \left\Vert z_{n}\right\Vert
_{L^{r}(B_{R}^{c})}^{r-1}+\left\Vert z_{n}\right\Vert
_{L^{q}(B_{R}^{c})}^{q-1}\right) \varepsilon .  \label{media}
\end{equation}

Since $z_{n}$ is bounded in $H^{1}\left( \mathbb{R}^{3}\right) ,$ the
sequences $\left\Vert z_{n}\right\Vert _{L^{r}(B_{R}^{c})}^{r-1}$ and $%
\left\Vert z_{n}\right\Vert _{L^{q}(B_{R}^{c})}^{q-1}$are bounded. Then, by (%
\ref{pepe}) and (\ref{media}), we easily get%
\begin{equation}
\underset{n\rightarrow \infty }{\lim }\left\vert K\left( u+w_{n}\right)
-K\left( u\right) -K\left( w_{n}\right) \right\vert \leq \varepsilon +M\cdot
\varepsilon .  \label{inno}
\end{equation}%
where $M$ is a suitable constant.

Since $\varepsilon $ is arbitrary, from (\ref{inno}) we get 
\begin{equation*}
\underset{n\rightarrow \infty }{\lim }\left\vert K\left( u+w_{n}\right)
-K\left( u\right) -K\left( w_{n}\right) \right\vert =0.
\end{equation*}

$\square $

\bigskip

In order to prove the coecitivity properties we need the following lemma:

\begin{lemma}
\label{coercive}Let the assumptions of Theorem \ref{exnsm} be satisfied.
Then $E\ $and $C$ defined by (\ref{rec}) and (\ref{enss}) satisfy the
coercivity assumption (EC-3).
\end{lemma}

\textbf{Proof.} By the Gagliardo-Nirenberg interpolation inequalities (see
e.g.\cite{niren}) there exists $b>0$ such that for any $u\in H^{1}(\mathbb{R}%
^{3})$ 
\begin{equation}
||u||_{L^{p}}^{p}\leq b||u||_{L^{2}}^{r}||\nabla u||_{L^{2}}^{q}
\label{main}
\end{equation}

where $q=3p\left( \frac{1}{2}-\frac{1}{p}\right) $ and $r=p-q.$ By (\ref{W0}%
) $2<p<\frac{10}{3},$ then $q<2$ and $r>0$.

Then by H\"{o}lder inequality we have for $M>0$%
\begin{align*}
||u||_{L^{p}}^{p}& \leq bM||u||_{L^{2}}^{r}\frac{1}{M}||\nabla
u||_{L^{2}}^{q} \\
& \leq \frac{1}{\gamma ^{\prime }}\left( bM||u||_{L^{2}}^{r}\right) ^{\gamma
^{\prime }}+\frac{1}{\gamma }\left( \frac{1}{M}||\nabla
u||_{L^{2}}^{q}\right) ^{\gamma } \\
& =\frac{\left( b_{p}M\right) ^{\gamma ^{\prime }}}{\gamma ^{\prime }}%
||u||_{L^{2}}^{r\gamma ^{\prime }}+\frac{1}{\gamma M^{\gamma }}||\nabla
u||_{L^{2}}^{q\gamma }.
\end{align*}%
Now chose $\gamma =\frac{2}{q}$ and $M=$ $\left( \frac{2c}{\gamma }\right)
^{1/\gamma }$, where $c$ is the constant in assumption (\ref{W0}), so that%
\begin{equation*}
||u||_{L^{p}}^{p}\leq \frac{\left( bM\right) ^{\gamma ^{\prime }}}{\gamma
^{\prime }}||u||_{L^{2}}^{r\gamma ^{\prime }}+\frac{1}{2c}||\nabla
u||_{L^{2}}^{2}.
\end{equation*}%
Then 
\begin{equation}
c||u||_{L^{p}}^{p}\leq a||u||_{L^{2}}^{2s}+\frac{1}{2}||\nabla
u||_{L^{2}}^{2}  \label{none}
\end{equation}%
where%
\begin{equation*}
a=\frac{c\left( bM\right) ^{\gamma ^{\prime }}}{\gamma ^{\prime }};\ \ \ s=%
\frac{r\gamma ^{\prime }}{2}.
\end{equation*}%
So using (\ref{W0}) , (\ref{none}) and setting 
\begin{equation}
\mathbf{F}^{2}=\Theta ^{2}+\mathbf{E}^{2}+\mathbf{H}^{2},  \label{somma}
\end{equation}%
we have for any $\mathbf{u}=\left( u,\Theta ,\mathbf{E},\mathbf{H}\right)
\in X$ 
\begin{align}
E(\mathbf{u})+aC(\mathbf{u})^{s}& =\frac{1}{2}\left\Vert \mathbf{u}%
\right\Vert ^{2}+\int V(x)u^{2}+\int N(u)+a||u||_{L^{2}}^{2s} \\
& \geq \frac{1}{2}||\nabla u||_{L^{2}}^{2}+\int \left( E_{0}u^{2}+\frac{%
\mathbf{F}^{2}}{2}\right) +\int N(u)+a||u||_{L^{2}}^{2s}  \notag \\
& \geq \frac{1}{2}||\nabla u||_{L^{2}}^{2}+\int \left( E_{0}u^{2}+\frac{%
\mathbf{F}^{2}}{2}\right) -c\int \left\vert u\right\vert
^{p}+a||u||_{L^{2}}^{2s}  \label{new} \\
& \geq \int \left( E_{0}u^{2}+\frac{\mathbf{F}^{2}}{2}\right) .
\label{mibis}
\end{align}%
Observe that, since $p>2,$ we have $s>1.$ So (EC-3)(i) is satisfied. Now we
prove that also (EC-3)(ii) is satisfied.

Let $\mathbf{u}_{n}$ $=\left( u_{n},\Theta _{n},\mathbf{E}_{n},\mathbf{H}%
_{n}\right) \in X$ be a sequence such that%
\begin{equation}
\left\Vert \mathbf{u}_{n}\right\Vert ^{2}=\int |\nabla u_{n}|^{2}+\int
\left( 2E_{0}u^{2}+\mathbf{F}_{n}^{2}\right) \rightarrow \infty .
\label{coco}
\end{equation}%
Now distinguish two cases:

- Assume first that $\int \left( 2E_{0}u^{2}+\mathbf{F}_{n}^{2}\right) $ is
unbounded. Then by (\ref{mibis}), we have (up to a subsequence) 
\begin{equation*}
E(\mathbf{u}_{n})+aC(\mathbf{u}_{n})^{s}\rightarrow \infty
\end{equation*}%
So in this case (EC-3$)$ (ii) is satisfied.

- Assume now that $\int \left( 2E_{0}u_{n}^{2}+\mathbf{F}_{n}^{2}\right) $
is bounded and set 
\begin{equation*}
d=\sup \left\Vert u_{n}\right\Vert _{L^{2}}^{p}.
\end{equation*}%
So by (\ref{main}) we have 
\begin{equation}
\left\Vert u_{n}\right\Vert _{L^{p}}^{p}\leq c_{1}||\nabla
u_{n}||_{L^{2}}^{q}\text{ where }c_{1}=bd.  \label{poi}
\end{equation}%
Since $\int \left( 2E_{0}u^{2}+\mathbf{F}_{n}^{2}\right) $ is bounded, by (%
\ref{coco}) we get 
\begin{equation}
\int |\nabla u_{n}|^{2}\rightarrow \infty .  \label{key}
\end{equation}%
On the other hand by (\ref{new}), we have 
\begin{equation*}
E(\mathbf{u}_{n})+aC(\mathbf{u}_{n})^{s}\geq \frac{1}{2}||\nabla
u_{n}||_{L^{2}}^{2}-c\int \left\vert u_{n}\right\vert ^{p}\geq (\text{by (%
\ref{poi}))}
\end{equation*}%
\begin{equation}
\frac{1}{2}||\nabla u_{n}||_{L^{2}}^{2}-c_{2}\text{ where }c_{2}=cc_{1}.
\label{koi}
\end{equation}%
Clearly (\ref{key}) and (\ref{koi}) prove that (EC-3)(ii) holds.

Now let us prove (EC-3)(iii). Let $\mathbf{u}_{n}$ $=\left( u_{n},\Theta
_{n},\mathbf{E}_{n},\mathbf{H}_{n}\right) \in X$ be a bounded sequence such
that $E(\mathbf{u}_{n})+aC(\mathbf{u}_{n})^{s}\rightarrow 0,$ then by (\ref%
{mibis}) we have%
\begin{equation}
\int \left( 2E_{0}u^{2}+\mathbf{F}_{n}^{2}\right) \rightarrow 0
\label{partone}
\end{equation}%
and hence%
\begin{equation}
\int u^{2}\rightarrow 0  \label{partina}
\end{equation}%
Then, in order to show that $\left\Vert \mathbf{u}_{n}\right\Vert $ $%
\rightarrow $ $0$ it remains to prove that%
\begin{equation}
||\nabla u_{n}||_{L^{2}}^{2}\rightarrow 0.  \label{norma}
\end{equation}%
Since $u_{n}$ is bounded in $H^{1}(\mathbb{R}^{3}),$ by (\ref{main}) and (%
\ref{partone}), we get 
\begin{equation}
\int \left\vert u_{n}\right\vert ^{p}\rightarrow 0.  \label{parte}
\end{equation}%
Since $E(\mathbf{u}_{n})+aC(\mathbf{u}_{n})^{s}\rightarrow 0$ and by
assumption \textbf{(}\ref{W0}\textbf{)}, we have%
\begin{eqnarray}
0 &=&\lim (E(\mathbf{u}_{n})+aC(\mathbf{u}_{n})^{s})  \label{comp} \\
&\geq &\lim \sup \left[ \frac{1}{2}||\nabla u_{n}||_{L^{2}}^{2}+E_{0}\int
\left\vert u_{n}\right\vert ^{2}-c\int \left\vert u_{n}\right\vert
^{p}+a||u_{n}||_{L^{2}}^{2s}\right] \\
&=&\lim \sup \left( \frac{1}{2}||\nabla u_{n}||_{L^{2}}^{2}+D_{n}\right)
\end{eqnarray}%
where%
\begin{equation}
D_{n}=E_{0}\int \left\vert u_{n}\right\vert ^{2}-c\int \left\vert
u_{n}\right\vert ^{p}+a||u_{n}||_{L^{2}}^{2s}.  \label{partissima}
\end{equation}%
By (\ref{partina}) and (\ref{parte}) $D_{n}\rightarrow 0.$ So by (\ref{comp}%
) we deduce (\ref{norma}).

$\square $

In the following we will verify that the hylomorphy condition (\ref{hh}) is
satisfied. For $\mathbf{u}=\left( u,\Theta ,\mathbf{E},\mathbf{H}\right) \in
X,$ we set

\begin{equation}
\left\Vert \mathbf{u}\right\Vert _{\sharp }=\left\Vert \left( u,\Theta ,%
\mathbf{E},\mathbf{H}\right) \right\Vert _{\sharp }=\left\Vert u\right\Vert
_{L^{t}},\text{ }2<t<6  \label{semina}
\end{equation}

\begin{equation*}
\Lambda _{0}:=\ \inf \left\{ \lim \inf \ \Lambda (\mathbf{u}_{n})\ |\ 
\mathbf{u}_{n}\ \text{is a vanishing sequence}\right\} ,\text{ }\Lambda
_{\sharp }=\ \underset{\left\Vert \mathbf{u}\right\Vert _{\sharp
}\rightarrow 0}{\lim \inf }\Lambda (\mathbf{u}).
\end{equation*}%
First of all we prove the following:

\begin{lemma}
\label{nonvanishing2}\textit{\textbf{\ }}The seminorm $\left\Vert \mathbf{u}%
\right\Vert _{\sharp }$ defined by (\ref{semina}) satisfies the property (%
\ref{seminorm}), namely, if $\mathbf{u}_{n}=\left( u_{n},\Theta _{n},\mathbf{%
E}_{n},\mathbf{H}_{n}\right) \ $is vanishing (see Definition \ref{na}), then 
$\left\Vert \mathbf{u}_{n}\right\Vert _{\sharp }=\left\Vert u_{n}\right\Vert
_{L^{t}}\rightarrow 0.$
\end{lemma}

\textbf{Proof. }For $j\in \mathbb{Z}^{3}$ we set 
\begin{equation*}
Q_{j}=A\left( j+Q^{0}\right) =\left\{ Aj+Aq:q\in Q^{0}\right\}
\end{equation*}%
where $Q^{0}$ is now the cube defined as follows 
\begin{equation*}
Q^{0}=\left\{ \left( x_{1},..,x_{n}\right) \in \mathbb{R}^{3}:0\leq
x_{i}<1\right\} \text{.}
\end{equation*}%
Now let $x\in \mathbb{R}^{3}$ and set $y=A^{-1}(x).$ Clearly there exist $%
q\in Q^{0}$ and $j\in \mathbb{Z}^{3}$ such that $y=j+q.$ So 
\begin{equation*}
x=Ay=A(j+q)\in Q_{j}.
\end{equation*}%
Then we conclude that 
\begin{equation*}
\mathbb{R}^{3}=\dbigcup\limits_{j}Q_{j}.
\end{equation*}%
Let $u_{n}\ $be a bounded sequence in $H^{1}\left( \mathbb{R}^{3}\right) $
such that, up to a subsequence, $\left\Vert u_{n}\right\Vert _{L^{t}}\geq
a>0.$ We need to show that $u_{n}$ is not vanishing. Then, if $L$ is the
constant for the Sobolev embedding $H^{1}\left( Q_{j}\right) \subset
L^{t}\left( Q_{j}\right) $ and $\left\Vert u_{n}\right\Vert _{H^{1}}^{2}\leq
M,$ we have 
\begin{align*}
0& <a^{t}\leq \int \left\vert u_{n}\right\vert
^{t}=\sum_{j}\int_{Q_{j}}\left\vert u_{n}\right\vert ^{t}=\sum_{j}\left\Vert
u_{n}\right\Vert _{L^{t}\left( Q_{j}\right) }^{t-2}\left\Vert
u_{n}\right\Vert _{L^{t}\left( Q_{j}\right) }^{2} \\
& \leq \ \left( \underset{j}{\sup }\left\Vert u_{n}\right\Vert _{L^{t}\left(
Q_{j}\right) }^{t-2}\right) \cdot \sum_{j}\left\Vert u_{n}\right\Vert
_{L^{t}\left( Q_{j}\right) }^{2} \\
& \leq \ L\left( \underset{j}{\sup }\left\Vert u_{n}\right\Vert
_{L^{t}\left( Q_{j}\right) }^{t-2}\right) \cdot \sum_{j}\left\Vert
u_{n}\right\Vert _{H^{1}\left( Q_{j}\right) }^{2} \\
& =L\left( \underset{j}{\sup }\left\Vert u_{n}\right\Vert _{L^{t}\left(
Q_{j}\right) }^{t-2}\right) \left\Vert u_{n}\right\Vert _{H^{1}}^{2}\leq
LM\left( \underset{j}{\sup }\left\Vert u_{n}\right\Vert _{L^{t}\left(
Q_{j}\right) }^{t-2}\right) .
\end{align*}%
Then%
\begin{equation*}
\left( \underset{j}{\sup }\left\Vert u_{n}\right\Vert _{L^{t}\left(
Q_{j}\right) }\right) \geq \left( \frac{a^{t}}{LM}\right) ^{1/(t-2)}
\end{equation*}%
Then, for any $n,$ there exists $j_{n}\in \mathbb{Z}^{3}$ such that 
\begin{equation}
\left\Vert u_{n}\right\Vert _{L^{t}\left( Q_{j_{n}}\right) }\geq \alpha >0.
\label{caca}
\end{equation}%
Then, if we set $Q=AQ^{0},$we easily have 
\begin{equation}
\left\Vert T_{j_{n}}u_{n}\right\Vert _{L^{t}(Q)}=\left\Vert u_{n}\right\Vert
_{L^{t}(Q_{j_{n}})}\geq \alpha >0.  \label{chicco}
\end{equation}

Since $u_{n}$ is bounded, also $T_{j_{n}}u_{n}$ is bounded in $H^{1}(\mathbb{%
R}^{3}).$ Then we have, up to a subsequence, that $T_{j_{n}}u_{n}%
\rightharpoonup u_{0}$ weakly in $H^{1}(\mathbb{R}^{3})$ and hence strongly
in $L^{t}(Q)$. By (\ref{chicco}), $u_{0}\neq 0.$

$\square $

\bigskip

By (\ref{ens}) and (\ref{somma}) the hylenic ratio takes the following form:%
\begin{equation}
\Lambda (u)=\frac{\int \left( \frac{1}{2}\left\vert \nabla u\right\vert
^{2}+E_{0}u^{2})+V(x)\left\vert u\right\vert ^{2}+\frac{\mathbf{F}^{2}}{2}%
\right) dx+\int N(u)}{\int \left\vert u\right\vert ^{2}dx}  \label{formino}
\end{equation}

\begin{lemma}
\label{preparatorio}If the assumptions of Theorem \ref{exnsm} are satisfied,
then for $2<t<6,$ we have%
\begin{equation*}
\underset{u\in H^{1},\left\Vert u\right\Vert _{L^{t}}\rightarrow 0}{\lim
\inf }\Lambda (u)\geq E_{0}
\end{equation*}
\end{lemma}

\textbf{Proof. }Clearly by (\ref{formino})

\begin{equation*}
\underset{u\in H^{1},\left\Vert u\right\Vert _{L^{t}}\rightarrow 0}{\lim
\inf }\Lambda (u)=\ \underset{u\in H^{1},\left\Vert u\right\Vert
_{L^{t}}=1,\varepsilon \rightarrow 0}{\lim \inf }\Lambda (\varepsilon u)
\end{equation*}%
\begin{equation*}
\geq E_{0}+\underset{u\in H^{1},\left\Vert u\right\Vert
_{L^{t}}=1,\varepsilon \rightarrow 0}{\lim \inf }\frac{\int N(\varepsilon
\psi )}{\varepsilon ^{2}\int \left\vert u\right\vert ^{2}}.
\end{equation*}%
So the proof of Lemma will be achieved if we show that%
\begin{equation}
\underset{u\in H^{1},\left\Vert u\right\Vert _{L^{t}}=1,\varepsilon
\rightarrow 0}{\lim \inf }\frac{\int N(\varepsilon u)}{\varepsilon ^{2}\int
\left\vert u\right\vert ^{2}}=0.  \label{resto}
\end{equation}%
By (\ref{Wp}) and ( \ref{W0}) we have%
\begin{equation}
-cs^{p}\leq N(s)\leq \bar{c}(s^{q}+s^{r})  \label{zerobis}
\end{equation}%
where $c,\bar{c}$ are positive constants and $q,r$ belonging to the interval 
$(2,2^{\ast }).$ Then by (\ref{zerobis}) we have%
\begin{equation}
-cA\varepsilon ^{p-2}\leq \underset{\left\Vert u\right\Vert _{L^{t}}=1}{\inf 
}\frac{\int N(\varepsilon u)}{\varepsilon ^{2}\int \left\vert u\right\vert
^{2}}\leq \bar{c}B(\varepsilon ^{q-2}+\varepsilon ^{r-2})  \label{uno}
\end{equation}%
where%
\begin{equation*}
A=\underset{u\in H^{1}\text{ }\left\Vert u\right\Vert _{L^{t}}=1}{\inf }%
\frac{\int \left\vert u\right\vert ^{p}}{\int \left\vert u\right\vert ^{2}},%
\text{ }B=\underset{u\in H^{1}\text{ }\left\Vert u\right\Vert _{L^{t}}=1}{%
\inf }\frac{\int \left( \left\vert u\right\vert ^{q}+\left\vert u\right\vert
^{r}\right) }{\int \left\vert u\right\vert ^{2}}.
\end{equation*}%
By (\ref{uno}) we easily get (\ref{resto}).

$\square $\bigskip

Now we can give an estimate of $\Lambda _{0}$ (see (\ref{hylo})).

\begin{corollary}
\label{interm} If the assumptions of Theorem \ref{exnsm} are satisfied, then%
\begin{equation*}
E_{0}\leq \Lambda _{0}.
\end{equation*}%
\bigskip\ 
\end{corollary}

\textbf{Proof.} By Proposition \ref{diesis}, Lemma \ref{nonvanishing2} and
Lemma \ref{preparatorio}

\begin{equation*}
\Lambda _{0}\geq \ \underset{\left\Vert u\right\Vert _{L^{t}}\rightarrow 0}{%
\Lambda _{\sharp }=\lim \inf }\ \Lambda (u)\geq E_{0}.
\end{equation*}

$\square $

\bigskip

\begin{lemma}
\label{buono}Let $W$ and $V$ satisfy assumptions (\ref{imis}), (\ref{W1}), (%
\ref{Wp}), (\ref{W0}), (\ref{V0}), (\ref{V1'}). Then, if $q$ is sufficiently
small, the hylomorphy condition (\ref{hh}) holds, namely%
\begin{equation}
\underset{\mathbf{u}\in X}{\inf }\Lambda (\mathbf{u})<\Lambda _{0}.
\label{conca}
\end{equation}
\end{lemma}

\textbf{Proof. }Clearly, by corollary \ref{interm}, in order to prove (\ref%
{conca}) it will be enough to show that for $q$ sufficiently small we have 
\begin{equation}
\underset{\mathbf{u}\in X}{\inf }\Lambda (\mathbf{u})<E_{0}.  \label{cota}
\end{equation}%
Taking $q$ sufficiently small, we will be able to construct $\mathbf{u}\in X$
such that $\Lambda (\mathbf{u})<E_{0}.$

Let $R>0$ and take $\mathbf{u}_{R}\mathbf{=}\left( u_{R},0,\mathbf{-}\nabla
\varphi _{R},0\right) ,$ where $u_{R}$ is defined by

\begin{equation*}
u_{R}=\left\{ 
\begin{array}{cc}
s_{0} & if\;\;|x|<R \\ 
0 & if\;\;|x|>R+1 \\ 
\frac{|x|}{R}s_{0}-(\left\vert x\right\vert -R)\frac{R+1}{R}s_{0} & 
if\;\;R<|x|<R+1%
\end{array}%
.\right.
\end{equation*}%
and $\varphi _{R}$ solves the equation 
\begin{equation}
-\Delta \varphi _{R}=qu_{R}^{2}.  \label{pois}
\end{equation}%
Take $\mathbf{u}_{R}\mathbf{=}\left( u_{R},0,\mathbf{-}\nabla \varphi
_{R},0\right) .$ Clearly, by definition of $X$\textbf{\ }(see(\ref{face}))
we have $\mathbf{u}_{R}\in X.$ Then 
\begin{equation*}
\int \left\vert \nabla u_{R}\right\vert ^{2}dx=O(R^{2}),\int \left\vert
u_{R}\right\vert ^{2}dx=O(R^{3}),
\end{equation*}%
so that%
\begin{equation}
\frac{\int \left[ \frac{1}{2}\left\vert \nabla u_{R}\right\vert ^{2}+\left(
E_{0}+V\right) u_{R}^{2}\right] dx}{\int u_{R}^{2}}\leq E_{0}+V_{0}+O\left( 
\frac{1}{R}\right) .  \label{pinco}
\end{equation}%
Moreover 
\begin{equation*}
\int N(u_{R})dx=N(s_{0})m(B_{R})+\int_{B_{R+1}\backslash B_{R}}N(u_{R}).
\end{equation*}%
where $m(A)$ denotes the measure of $A.$ So%
\begin{eqnarray}
\frac{\int N(u_{R})dx}{\int u_{R}^{2}} &\leq &\frac{%
N(s_{0})m(B_{R})+c_{1}R^{2}}{\int u_{R}^{2}}\leq (\text{ since }N(s_{0})<0)
\label{palle} \\
&\leq &\frac{N(s_{0})m(B_{R})}{s_{0}^{2}m(B_{R+1})}+\frac{c_{1}R^{2}}{%
s_{0}^{2}m(B_{R})}=\frac{N(s_{0})}{s_{0}^{2}}\left( \frac{R}{R+1}\right)
^{3}+\frac{c_{2}}{R}.  \notag
\end{eqnarray}%
Then, since $\Theta =\mathbf{H}=0,$ by (\ref{formino}), (\ref{pinco}) and (%
\ref{palle}) we get%
\begin{eqnarray}
&&\Lambda (\mathbf{u}_{R}) \\
&=&\ \frac{\int \left( \frac{1}{2}\left\vert \nabla u_{R}\right\vert
^{2}+\left( E_{0}+V(x)\right) u_{R}^{2}\right) dx}{\int u_{R}^{2}dx}+\frac{%
\int N(u_{R})dx}{\int u_{R}^{2}dx}+\frac{\frac{1}{2}\int \left\vert \nabla
\varphi _{R}\right\vert ^{2}}{\int u_{R}^{2}}\leq  \notag \\
&\leq &E_{0}+V_{0}+\frac{N(s_{0})}{s_{0}^{2}}\left( \frac{R}{R+1}\right)
^{3}+\frac{c_{2}}{R}+\frac{\frac{1}{2}\int \left\vert \nabla \varphi
_{R}\right\vert ^{2}}{\int u_{R}^{2}}.  \label{pinko}
\end{eqnarray}%
Now we will estimate the term containing $\varphi _{R}$ in (\ref{pinko}).
Observe that $u_{R}^{2}$ has radial symmetry and that the electric field
outside any spherically symmetric charge distribution is the same as if all
of the charge were concentrated into a point. So $\left\vert \nabla \varphi
_{R}\left( r\right) \right\vert $ corresponds to the strength of an
electrostatic field at distance $r,$ created by an electric charge given by 
\begin{equation*}
\left\vert C_{el}\right\vert =\dint\limits_{\left\vert x\right\vert \leq
r}qu_{R}^{2}dx=4\pi \dint\limits_{0}^{r}qu_{R}^{2}v^{2}dv
\end{equation*}%
and located at the origin. So we have%
\begin{equation*}
\left\vert \nabla \varphi _{R}\left( r\right) \right\vert =\frac{\left\vert
C_{el}\right\vert }{r^{2}}\left\{ 
\begin{array}{cc}
=\frac{4}{3}\pi qs_{0}^{2}r & if\ r<R \\ 
\leq \frac{4}{3}\pi qs_{0}^{2}\frac{(R+1)^{3}}{r^{2}} & if\ r\geq R%
\end{array}%
\right. .
\end{equation*}%
Then%
\begin{align*}
\int \left\vert \nabla \varphi _{R}\right\vert ^{2}dx& \leq
c_{3}q^{2}s_{0}^{4}\left( \int_{r<R}r^{2}dr+\int_{r>R}\frac{(R+1)^{6}}{r^{2}}%
dr\right) \\
& \leq c_{4}q^{2}s_{0}^{4}\left( R^{3}+\frac{(R+1)^{6}}{R}\right) \leq
c_{5}q^{2}s_{0}^{4}R^{5}.
\end{align*}%
So%
\begin{equation}
\frac{\frac{1}{2}\int \left\vert \nabla \varphi _{R}\right\vert ^{2}}{\int
u_{R}^{2}}\leq \frac{c_{6}\int \left\vert \nabla \varphi _{R}\right\vert ^{2}%
}{s_{0}^{2}R^{3}\ }\leq c_{7}q^{2}s_{0}^{2}R^{2}.  \label{bina}
\end{equation}%
By (\ref{bina}) and (\ref{pinko}), we get%
\begin{equation}
\Lambda (\mathbf{u}_{R})\leq E_{0}+V_{0}+\frac{N(s_{0})}{s_{0}^{2}}\left( 
\frac{R}{R+1}\right) ^{3}+\frac{c_{2}}{R}+c_{7}q^{2}s_{0}^{2}R^{2}.
\label{tania}
\end{equation}%
Since by our assumptions 
\begin{equation*}
\frac{N(s_{0})}{s_{0}^{2}}<-V_{0}
\end{equation*}%
for $R$ large we get 
\begin{equation}
V_{0}+\frac{N(s_{0})}{s_{0}^{2}}\left( \frac{R}{R+1}\right) ^{3}+\frac{c_{2}%
}{R}<0  \label{tonio}
\end{equation}%
So, if $q$ is small enough, by (\ref{tania}) and (\ref{tonio}) we get

\begin{equation*}
\Lambda (\mathbf{u}_{R})<E_{0}.
\end{equation*}

$\square $

\bigskip

\textbf{Proof of Theorem} \ref{exnsm}. We shall show that all the
assumptions of theorem \ref{astra1} are satisfied. Assumptions (EC-0),
(EC-1) , are clearly satisfied. By Lemma \ref{split} and Lemma \ref{coercive}
also the splitting property (EC-2) and the coercivity property (EC-3) hold.
By Lemma \ref{buono} the hylomorphy condition (\ref{hh}) holds. Finally also
the assumption (\ref{poo}) is satisfied. In fact it is immediate to see that 
\begin{equation*}
E^{\prime }\left( u,\Theta ,\mathbf{E},\mathbf{H}\right) =0\Longrightarrow
\Theta =\mathbf{E}=\mathbf{H=}0
\end{equation*}%
\begin{equation*}
C^{\prime }(u,\Theta ,\mathbf{E},\mathbf{H)}=0\Longrightarrow C^{\prime }(u%
\mathbf{)}=0\Longrightarrow u=0.
\end{equation*}%
$\square $

\bigskip

\textbf{Proof of Theorem }\ref{nsmbis} Let $\mathbf{u}_{\delta }=\left(
u_{\delta },\Theta _{\delta },\mathbf{E}_{\delta },\mathbf{H}_{\delta
}\right) $ be an hylomorphic soliton for NSM. So there exists a constant $%
\sigma $ such $\mathbf{u}_{\delta }$ minimizes the energy $E$ (see (\ref{ens}%
)) on the manifold%
\begin{equation*}
\mathfrak{M}_{\sigma }=\left\{ \mathbf{u=}\left( u,\Theta ,\mathbf{E},%
\mathbf{H}\right) \in X:C(\mathbf{u)=}\int u^{2}dx=\sigma \right\}
\end{equation*}%
where 
\begin{equation*}
X=\left\{ \mathbf{u}=\left( u,\Theta ,\mathbf{E},\mathbf{H}\right) \in
H^{1}\left( \mathbb{R}^{3}\right) \times L^{2}\left( \mathbb{R}^{3}\right)
^{9}:\nabla \cdot \mathbf{E}=qu^{2}\right\} .
\end{equation*}%
Since $\mathbf{u}_{\delta }$ $=\left( u_{\delta },\Theta _{\delta },\mathbf{E%
}_{\delta },\mathbf{H}_{\delta }\right) $ minimizes the energy $E$ on $%
\mathfrak{M}_{\sigma },$ we have $\Theta _{\delta }=\mathbf{H}_{\delta }=0,$
then%
\begin{equation*}
\mathbf{u}_{\delta }=\left( u_{\delta },0,\mathbf{E}_{\delta }\mathbf{,}%
0\right)
\end{equation*}%
If we set $\mathbf{E}=-\nabla \varphi $, the constraint $\nabla \cdot 
\mathbf{E}=qu^{2}$ becomes 
\begin{equation}
-\Delta \varphi =qu^{2}.  \label{provvi}
\end{equation}%
So $\mathbf{u}_{\delta }$ is a critical point of $E$ on the manifold made up
by those $\mathbf{u}=$ $\left( u,0,-\nabla \varphi \mathbf{,}0\right) $
satisfying the constraints (\ref{provvi}) and%
\begin{equation}
C(\mathbf{u)=}\int u^{2}\ dx=\sigma .  \label{vic20}
\end{equation}%
Therefore, for suitable Lagrange multipliers $\omega \in \mathbb{R},$ $\xi
\in \mathcal{D}^{1,2}$ ($\mathcal{D}^{1,2}$ is the closure of $C_{0}^{\infty
}$ with respect to the norm $\left\Vert \nabla \varphi \right\Vert _{L^{2}}$%
), we have that $\mathbf{u}_{\delta }$ is a critical point of the free
functional 
\begin{eqnarray}
E_{\omega ,\xi }(\mathbf{u)} &=&E(\mathbf{u)}+\omega \left( \int
u^{2}-\sigma \right) +\left\langle \xi ,\Delta \varphi +qu^{2}\right\rangle 
\text{ }  \label{lager} \\
&=&\int \left( \frac{1}{2}\left\vert \nabla u\right\vert ^{2}+V(x)u^{2}+W(u)+%
\frac{1}{2}\left\vert \nabla \varphi \right\vert ^{2}\right) dx  \notag \\
&&+\omega \left( \int u^{2}-\sigma \right) +\left\langle \xi ,\Delta \varphi
+qu^{2}\right\rangle \text{ }  \notag
\end{eqnarray}%
where $\left\langle \ \cdot \ ,\ \cdot \ \right\rangle $ denotes the duality
map in $\mathcal{D}^{1,2}$ and $\mathbf{u}$ can be identified with $%
(u,\varphi )\in H^{1}\left( \mathbb{R}^{3}\right) \times $ $\mathcal{D}%
^{1,2}.$ $E_{\omega ,\xi }^{\prime }(\mathbf{u}_{\delta })=0$ gives the
equations%
\begin{eqnarray*}
\forall v &\in &H^{1}\left( \mathbb{R}^{3}\right) ,\ \ \left\langle \frac{%
\partial E_{\omega ,\xi }(\mathbf{u)}}{\partial u},v\right\rangle =0 \\
\forall \chi &\in &\mathcal{D}^{1,2}\left( \mathbb{R}^{3}\right) ,\ \
\left\langle \frac{\partial E_{\omega ,\xi }(\mathbf{u)}}{\partial \varphi }%
,\chi \right\rangle =0
\end{eqnarray*}%
namely%
\begin{equation*}
\forall v\in H^{1}\left( \mathbb{R}^{3}\right) ,\ \int \nabla u\cdot \nabla
v+2\left[ V(x)u+\frac{1}{2}W^{\prime }(u)+\omega u+q\xi u\right] v=0
\end{equation*}%
\begin{equation*}
\forall \chi \in \mathcal{D}^{1,2},\ \int \nabla \varphi \cdot \nabla \chi
+\left\langle \xi ,\Delta \chi \right\rangle =0.
\end{equation*}%
So, $u_{\delta },\varphi _{\delta }$ are weak solutions of the following
equations: 
\begin{align}
-\frac{1}{2}\Delta u_{\delta }+V(x)u_{\delta }+\frac{1}{2}W^{\prime
}(u_{\delta })+\omega u_{\delta }+q\xi u_{\delta }& =0  \label{svar1} \\
\Delta \varphi _{\delta }& =\Delta \xi .  \label{svar3}
\end{align}%
From (\ref{svar3}) we get $\xi =\varphi _{\delta },$ so (\ref{svar1}) becomes%
\begin{equation*}
-\Delta u_{\delta }+2V(x)u_{\delta }+W^{\prime }(u_{\delta })+2\omega
u_{\delta }+2q\varphi _{\delta }u_{\delta }=0.
\end{equation*}%
This equation and the constraint (\ref{provvi}) give the system (\ref{din})
and (\ref{dan}).

$\square $

\end{document}